\documentclass[11pt]{article}
\usepackage{amsmath}
\usepackage{anysize}

\marginsize{0.65in}{0.65in}{1in}{1in}
\usepackage{setspace}
\singlespace
\usepackage[dvips]{color}
\definecolor{myred}{rgb}{1,0.2,0}

\usepackage{stmaryrd}
\usepackage{amsfonts,amsmath,amssymb}
\usepackage{times}
\usepackage{color}
\usepackage[table]{xcolor}

\usepackage{multirow}
\usepackage{rotating}
\usepackage{bbm}
\usepackage{latexsym}

\usepackage{cite}

\usepackage{algorithm}
\usepackage{algorithmic}
\usepackage{booktabs}
\usepackage{url}
\usepackage{pifont}

\usepackage{psfrag}

\usepackage{colortbl}

\definecolor{shadingcolor}{rgb} {.87,    0.92,    .98}
\newcommand{\shadingbox}[1]{   %% shading for table
    \fboxsep 0pt
    \colorbox{shadingcolor}{
        {\hskip -2pt #1}\hskip -2.5pt
    }
}

\newcommand{\minitab}[2][l]{\begin{tabular}{@{}#1}#2\end{tabular}}

\newcommand{\tensor}[1]{\boldsymbol{\mathscr{\MakeUppercase{#1}}}} %tensor
\newcommand{\tA}{\tensor{A}}
\newcommand{\tB}{\tensor{B}}
\newcommand{\tC}{\tensor{C}}
\newcommand{\tD}{\tensor{D}}

\newcommand{\tG}{\tensor{G}}
\newcommand{\tH}{\tensor{H}}

\newcommand{\tS}{\tensor{S}}

\newcommand{\tX}{\tensor{X}}
\newcommand{\tY}{\tensor{Y}}

\input{def1.set}

\usepackage{soul}
\usepackage{xcolor}

\setstcolor{red}

\newcommand{\remove}[1]{}
\newcommand{\add}[1]{{{{#1}}}}

\newcommand{\replace}[2]{\remove{#1}  \add{#2} }

\usepackage{mathpazo}

\usepackage[mathscr]{eucal}

\graphicspath{{BCDCPDvsICA/}}

\newcommand{\matn}[2][n]{\ensuremath{\mathbf{#2}^{(#1)}}}
\renewcommand{\figurename}{Figure~}

\usepackage[font=small,labelfont=bf,tableposition=top]{caption}
\usepackage[font=footnotesize]{subcaption}

\usepackage{rotating}

\date{}
\begin{document}
%

%\title{Multi-Way Component Analysis: A Tensor Signal Processing Perspective}
\title{Tensor Decompositions for Signal Processing Applications \\
{\Large From Two-way to Multiway Component Analysis\footnote{Preprint of feature paper which will be published in the IEEE Signal Processing Magazine (2014).}
}}

\author{A. Cichocki, D. Mandic,  A-H. Phan, C. Caiafa,  G. Zhou, Q. Zhao,  and\\ L. De Lathauwer}

\maketitle

{\bf Summary} 
 The widespread use of multi-sensor technology and the emergence of big datasets has highlighted the limitations of standard flat-view matrix models and the necessity to move
 towards  more versatile data analysis tools. We show that higher-order tensors (i.e., multiway arrays) enable such a fundamental paradigm shift towards models that are essentially polynomial and whose uniqueness, unlike the matrix methods, is guaranteed under very mild and natural conditions.
 Benefiting from the power of multilinear algebra as their mathematical backbone, data analysis techniques using tensor decompositions are shown to have great flexibility in the choice of constraints that match data properties, and to find more general latent components in the data than matrix-based methods.
 % and help uniqueness.
 A comprehensive introduction to tensor decompositions is provided from a signal processing perspective, starting from the algebraic foundations, via basic  Canonical Polyadic and Tucker models, through  to advanced cause-effect  and multi-view data analysis schemes. We show that tensor decompositions enable natural generalizations of some commonly used signal processing paradigms, such as canonical correlation and subspace techniques, signal separation, linear regression, feature extraction and classification. We also cover computational aspects, and point out how ideas from compressed sensing and scientific computing may be used for addressing the otherwise unmanageable storage and manipulation problems associated with big datasets. The concepts are supported by illustrative real world case studies illuminating the \replace{blessings}{benefits} of the tensor framework, as efficient and promising tools for modern signal processing, data analysis and machine learning applications; these \replace{blessings}{benefits} also extend to vector/matrix data through tensorization.


\section{Introduction}

{\bf Historical notes.} The roots of multiway analysis can be traced back to studies
of homogeneous polynomials in the 19th century, contributors include
Gauss, Kronecker, Cayley, Weyl and Hilbert ---  in modern day
interpretation these are  fully symmetric tensors.
Decompositions of non-symmetric tensors have been studied since the early 20th century \cite{Hitchcock1927}, whereas
the benefits of using more than two matrices in factor analysis\cite{Cattell} became apparent in several communities since the 1960s. The Tucker decomposition
for tensors was introduced in psychometrics
\cite{tucker64extension,Tucker1966}, while the Canonical Polyadic Decomposition (CPD)
was independently rediscovered and put into an application context under the names of
Canonical Decomposition (CANDECOMP) in psychometrics \cite{PARAFAC1970Carroll} and Parallel Factor Model (PARAFAC) in linguistics
\cite{Harshman}.
Tensors were subsequently adopted in  diverse branches of data analysis such as chemometrics, food industry and social sciences \cite{Smilde,Kroonenberg}.
When it comes to Signal Processing, the early 1990s saw a considerable interest in
Higher-Order Statistics (HOS) \cite{Nikias-Petropulu-HOS} and it was soon realized that for the multivariate case HOS are effectively higher-order tensors;
 indeed, algebraic approaches to Independent Component Analysis (ICA)  using  HOS \cite{cardoso1993blind,comon94,Comon-Jutten2010} were inherently tensor-based. Around 2000, it was realized that the Tucker decomposition represents a MultiLinear Singular Value Decomposition (MLSVD) \cite{Lathauwer01}.
Generalizing the matrix SVD, the workhorse of numerical
linear algebra, the MLSVD spurred the interest in tensors in applied mathematics and scientific computing in very high dimensions \cite{beylkin,hTucker,OseledetsTT11}. In parallel, CPD was successfully adopted as a  tool for sensor array processing  and deterministic signal separation in wireless communication  \cite{Sidiropoulos00Bro,SidiropoulosGianna}. Subsequently, tensors have been used in audio, image and video processing, machine learning and biomedical applications, to name but a few.
The significant interest in tensors and their fast emerging applications are reflected in books \cite{Smilde,Kroonenberg,NMF-book,Comon-Jutten2010,Landsberg,Hackbush2012} and tutorial papers \cite{Acar2008,Kolda08,Comon-ALS09,Lu-2011,Morup11,Khoromskij-TT,Grasedyck-rev,ComonSPM2014} covering various aspects of multiway analysis.

{\bf From a matrix to a tensor.} Approaches to two-way (matrix) component analysis are well established, and include  Principal Component Analysis (PCA), Independent Component Analysis (ICA), Nonnegative Matrix Factorization (NMF) and Sparse Component Analysis (SCA) \cite{bruckstein2009sparse,NMF-book,Comon-Jutten2010}.
These techniques have become standard tools for e.g., blind source separation (BSS), feature extraction, or  classification.
On the other hand, large classes of data arising
from modern heterogeneous sensor modalities have a multiway character and are therefore
naturally represented by multiway arrays or tensors (see Section {\bf Tensorization}).

Early multiway data analysis approaches reformatted the data tensor as a matrix and resorted to methods developed for classical two-way analysis. However, such a ``flattened'' view of the world and the rigid assumptions inherent in two-way analysis are not always a good match for multiway data.
It is only through higher-order tensor decomposition that we have the opportunity to develop sophisticated models capturing multiple interactions and couplings, instead of standard  pairwise interactions.
In other words, we can only discover hidden components within multiway data if the analysis tools account for intrinsic multi-dimensional patterns present --- motivating the development of multilinear techniques.

In this article, we emphasize that tensor decompositions are not just matrix factorizations with additional subscripts --- multilinear algebra is much structurally richer than linear algebra. For example, even basic notions such as rank have a more subtle meaning,
uniqueness conditions of higher-order tensor decompositions are more relaxed and accommodating than those for  matrices \cite{Kruskal1977,Domanov12}, while matrices and tensors also have completely different geometric properties \cite{Landsberg}. This boils down to matrices representing linear transformations and quadratic forms, while tensors are connected with multilinear mappings and multivariate polynomials \cite{ComonSPM2014}.

\section{Notations and Conventions}

A tensor can be thought of as a multi-index numerical array, whereby the order of a tensor is the number of its ``{\emph{modes}}'' or ``{\emph{dimensions}}'', these may include space, time, frequency,  trials, classes, and dictionaries. A real-valued tensor of order $N$ is denoted by $\tA \in \Real^{I_1 \times I_2 \times \cdots \times I_N}$ and its entries by $a_{i_1,i_2,\ldots,i_N}$. Then, an $N \times 1$ vector $\ba$ is considered a tensor of order one, and an $N \times M$ matrix $\bA$ a tensor of order two. Subtensors are parts of the original data tensor, created when only a fixed subset of indices is used. Vector-valued subtensors are called {\it fibers}, defined  by fixing every index but one, and matrix-valued subtensors are called {\it  slices}, obtained by fixing all but two indices (see Table \ref{table_notation1}). Manipulation of tensors often requires their reformatting ({\emph{reshaping}});
%\remove{of tensors},
%\replace{especially we often need to reshape tensors to matrices; such process}{
\replace{the}{a} particular case of reshaping tensors to matrices is termed matrix unfolding or matricization (see Figure \ref{Fig:MWCA} (left)).
Note that a mode-$n$ multiplication of a tensor $\tA$ with a matrix  $\bB$ amounts to the multiplication of all mode-$n$ vector fibers with $\bB$, and that in linear algebra the tensor (or outer) product appears in the expression for a rank-1 matrix: $\ba \bb^T = \ba \circ \bb$.
Basic tensor notations are summarized in Table \ref{table_notation1}, while Table \ref{table_notation2} outlines several types of products used in this paper.

\minrowclearance 2ex
\begin{table*}[t!]
\caption{Basic notation.} \centering
 {\small  \shadingbox{
    \begin{tabular*}{\linewidth}[t]{@{\extracolsep{\fill}}ll} \hline
{$\tA, \; \bA, \; \ba, \; a$}& {tensor, matrix, vector, scalar} \\
$\bA = [\ba_1,\ba_2,\ldots,\ba_R]$ & {matrix  $\bA$ with column vectors $\ba_r$}  \\
%[1ex]
%$\bD = \mbox {diag} (\lambda_1,\lambda_2,\ldots,\lambda_R)$ & {diagonal matrix with $d_{rr}=\lambda_r$}  \\
$\ba(:,i_2,i_3,\ldots,i_N)$ & \minitab[p{.6\linewidth}]{fiber of tensor $\tA$ obtained by fixing all but one index}\\[-1ex]
$\bA(:,:,i_3,\ldots,i_N)$ & \minitab[p{.6\linewidth}]{matrix slice of tensor $\tA$ obtained by fixing all but two indices}\\[-1ex]
$\tA(:,:,:,i_4,\ldots,i_N)$ & \minitab[p{.6\linewidth}]{tensor slice of $\tA$ obtained by fixing some indices}\\[-1ex]
%CesarAddedOct17
$\tA(\mathcal{I}_1,\mathcal{I}_2,\ldots,\mathcal{I}_N)$ & \minitab[p{.6\linewidth}]{subtensor of $\tA$ obtained by restricting indices to belong to subsets $\mathcal{I}_n \subseteq \{ 1,2,\dots, I_n \}$}\\[-1ex]
$\bA_{(n)} \in \Real^{I_n \times I_1 I_2 \cdots I_{n-1} I_{n+1} \cdots I_N}$ & \minitab[p{.6\linewidth}]{mode-$n$ matricization of tensor $\tA \in \Real^{I_1 \times I_2 \times \cdots \times I_N}$ whose entry at row $i_n$ and column $(i_1 -1) I_2 \cdots I_{n-1} I_{n+1} \cdots I_N +  \cdots + (i_{N-1} -1)I_N + i_N$ is equal to $a_{i_1 i_2 \ldots i_N}$}
\\[-3ex]
$\vtr{\tA}\in \Real^{I_N I_{N-1} \cdots I_1}$ & \minitab[p{.6\linewidth}]{vectorization of tensor $\tA \in \Real^{I_1 \times I_2 \times \cdots \times I_N}$ with the entry at position $i_1 +\sum_{k=2}^N[(i_k-1)I_1 I_2 \cdots I_{k-1}]$ equal to $a_{i_1 i_2 \ldots i_N}$
%$(i_N -1) I_{N-1} I_{N-2} \cdots I_1 +  \cdots + (i_2 -1)I_1 + i_1$ equal to $a_{i_1 i_2 \ldots i_N}$
%with an entry at position $i_n$ and column $(i_N -1) I_{N} \cdots I_1 +  \cdots + (i_2 -1)I_1 + i_1$ equal to $a_{i_1 i_2 \ldots i_N}$}
}\\[1ex]
$\bD = \mbox {diag} (\lambda_1,\lambda_2,\ldots,\lambda_R)$ & {diagonal matrix with $d_{rr}=\lambda_r$}  \\
%$\bD = \diag_2(\lambda_1, \lambda_2, \ldots, \lambda_R)$ & diagonal matrix \\
$\tD = \diag_N(\lambda_1, \lambda_2, \ldots, \lambda_R)$ & diagonal tensor of order $N$ with $d_{rr \cdots r} = \lambda_r$ \\
${\bA}^{T}$,   ${\bA}^{-1}$, ${\bA}^{\dag}$ & {transpose, inverse, and Moore-Penrose pseudo-inverse} \\
\hline
    \end{tabular*}
   }}
%    \end{center}
\label{table_notation1}
\end{table*}
\minrowclearance 0ex

\minrowclearance 2ex
\begin{table*}[ht]
\caption{Definition of products.} %\centering
  {\small \shadingbox{
    \begin{tabular*}{\linewidth}[t]{@{\extracolsep{\fill}}ll} \hline
\\[-4em]
$\tC = \tA \times_n \bB$ & \minitab[p{.65\linewidth}]{\\[-3.2ex] mode-$n$ product of $\tA \in \Real^{I_1 \times I_2 \times \cdots \times I_N}$ and $\bB \in \Real^{J_n \times I_n}$ yields $\tC \in \Real^{I_1 \times \cdots \times I_{n-1} \times J_n \times I_{n+1} \times \cdots \times I_N}$ with entries $c_{i_1 \cdots i_{n-1} \, j_n \, i_{n+1} \cdots i_N} = \sum_{i_n=1}^{I_n} a_{i_1 \cdots i_{n-1} \, i_n \, i_{n+1} \cdots i_N} b_{j_n \, i_n}$ and matrix representation $\bC_{(n)} = \bB \, \bA_{(n)}$} \\
{$\tC = \llbracket \tA; \bB^{(1)}, \bB^{(2)}, \ldots, \bB^{(N)}\rrbracket$ } & {full multilinear product, $\tC =   {\tA} \times_1 \bB^{(1)} \times_2 \bB^{(2)}
\cdots \times_N \bB^{(N)}$ } \\ [-2ex]
$\tC = \tA \circ \tB$ & \minitab[p{.65\linewidth}]{\\[-3.2ex] tensor or outer product of $\tA \in \Real^{I_1 \times I_2 \times \cdots \times I_N}$ and $\tB \in \Real^{J_1 \times J_2 \times \cdots \times J_M}$ yields $\tC \in \Real^{I_1 \times I_2 \times \cdots \times I_N \times J_1 \times J_2 \times \cdots \times J_M}$ with entries $c_{i_1 i_2 \cdots i_N \,j_1 j_2 \cdots j_N} = a_{i_1 i_2 \cdots i_N} b_{j_1 j_2 \cdots j_M}$} \\ [-2ex]
$\tX = \ba^{(1)} \circ \ba^{(2)} \circ \cdots \circ \ba^{(N)}$ &\minitab[p{.7\linewidth}]{\\[-3.2ex] tensor or outer product of vectors $\ba^{(n)} \in \Real^{I_n}\;$ $(n=1,\ldots,N)$ yields a rank-1 tensor $\tX \in \Real^{I_1 \times I_2 \times \cdots \times I_N}$ with entries $x_{i_1 i_2 \ldots i_N} = a^{(1)}_{i_1} a^{(2)}_{i_2} \ldots a^{(N)}_{i_N} $}  \\ [-2ex]
{$\bC = \bA\otimes\bB$}&\minitab[p{.65\linewidth}]{Kronecker product of $\bA \in \Real^{I_1 \times I_2}$ and $\bB \in \Real^{J_1 \times J_2}$ yields $\bC \in \Real^{I_1 J_1 \times I_2 J_2}$ with entries $c_{(i_1-1) J_1 +j_1,(i_2-1) J_2 + j_2} = a_{i_1 i_2} \: b_{j_1 j_2}$} \\[-2ex]
{$\bC = \bA\odot \bB$}&\minitab[p{.65 \linewidth}]{Khatri-Rao product of $\bA=[\ba_1,\ldots,\ba_R] \in \Real^{I\times R}$ and $\bB=[\bb_1,\ldots,\bb_R] \in \Real^{J\times R}$ yields $\bC \in \Real^{I J\times R}$  with columns $\bc_r = \ba_r \otimes \bb_r$}\\ \hline
    \end{tabular*}
   }}
%    \end{center}
\label{table_notation2}
\end{table*}
\minrowclearance 0ex

\section{Interpretable Components in Two-Way Data Analysis}

The aim of blind source separation (BSS), factor analysis (FA) and latent variable analysis (LVA) is to decompose a data matrix $\bX \in \Real^{I \times J}$ into the factor matrices $\bA=[\ba_1,\ba_2,\ldots,\ba_R] \in \Real^{I \times R}$ and $\bB=[\bb_1,\bb_2,\ldots,\bb_R] \in \Real^{J \times R}$ as:
\be
 \label{2CA}
 \bX  &=& \bA \, \bD\, \bB^T + \bE =\sum_{r=1}^R \lambda_r \, \ba_r  \bb_r^{T} +\bE \nonumber \\
 &=& \sum_{r=1}^R \lambda_r \, \ba_r  \circ \bb_r +\bE,
\ee
where $ \bD = \mbox{\rm diag}(\lambda_1, \lambda_2, \ldots, \lambda_R)$ is a scaling (normalizing)  matrix, the columns of $\bB$ represent the unknown source signals (factors or latent variables depending on the tasks in hand), the columns of $\bA $ represent the associated mixing vectors (or factor loadings), while $\bE$ is noise due to an unmodelled data part  or model error. In other words, model (\ref{2CA}) assumes that the data matrix $\bX$ comprises hidden components $\bb_r$ ($r=1,2,\ldots,R$) that are mixed together in an unknown manner through coefficients $\bA$, or, equivalently, that data contain {\emph{factors}} that have an associated {\emph{loading}} for every data channel. Figure \ref{fig:CPD} (top) depicts the model (\ref{2CA}) as a dyadic decomposition, whereby the terms $\ba_r \circ \bb_r=\ba_r \bb_r^T$  are rank-1 matrices.

The well-known indeterminacies intrinsic to this model are: (i) arbitrary scaling %and counter scaling
of components, and (ii) permutation of the rank-1 terms. Another indeterminacy is related to the physical meaning of the factors: if the model in (\ref{2CA}) is unconstrained, it admits infinitely many combinations of $\bA$ and $\bB$. Standard matrix factorizations in linear algebra, such as the QR-factorization, Eigenvalue Decomposition (EVD), and Singular Value Decomposition (SVD), are only special cases of (\ref{2CA}), and owe their uniqueness to hard and restrictive constraints such as triangularity and orthogonality.
%
%On the other hand, when %\remove{the factors in (\ref{2CA}) exhibit }
%certain properties \add{of the factors in (\ref{2CA})} \remove{, these} can be represented by appropriate constraints \replace{on}{within the} model (\ref{2CA}), making possible \remove{their} unique estimation or extraction \add{of such factors}.
On the other hand, certain properties of the factors in (\ref{2CA}) can be represented by appropriate constraints, making possible unique
estimation or extraction of such factors.
\replace{Such}{These} constraints \replace{are, for instance,}{include} statistical independence, sparsity, nonnegativity, exponential structure, uncorrelatedness, constant modulus, finite alphabet, smoothness and unimodality. Indeed, the first four properties form the basis of Independent Component Analysis (ICA) \cite{Comon-Jutten2010,CiAm02,HyvICA2013}, Sparse Component Analysis (SCA) \cite{bruckstein2009sparse}, Nonnegative Matrix Factorization (NMF) \cite{NMF-book}, and harmonic retrieval  \cite{elad2004shape}.

\section{Tensorization --- Blessing of Dimensionality} \label{sect:tensorization}
While one-way (vectors) and two-way (matrices) algebraic structures were respectively introduced as natural representations for segments of scalar measurements and measurements on a grid, tensors were initially used purely for the mathematical benefits they provide in data analysis; for instance, it seemed natural to stack together excitation-emission spectroscopy matrices  in chemometrics into a third-order tensor \cite{Smilde}.

\begin{figure}[t!]
\begin{center}
\includegraphics[width=10.2cm]{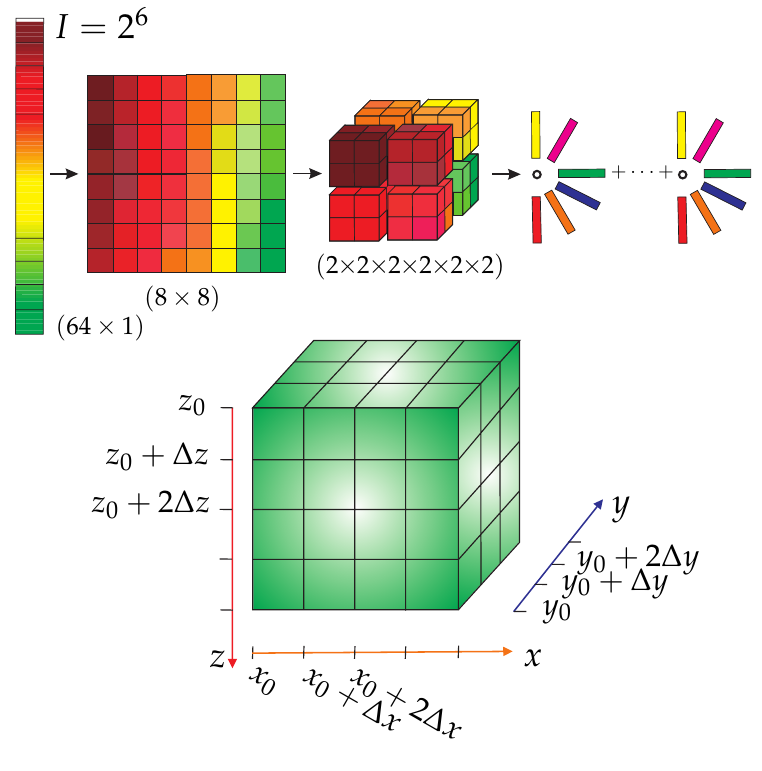}
\end{center}
\vspace{-0.3cm}
\caption{Construction of tensors. {\em Top:} Tensorization of a vector or matrix into the so-called {\em quantized format}; in scientific computing this facilitates super-compression of large-scale vectors or matrices. %(e.g., through a polyadic or a TT decomposition).
{\em Bottom:} Tensor formed through the discretization of a trivariate function $f(x,y,z)$.}
\label{Fig:Tensorization}
\end{figure}

%\begin{figure}[t!]
%\begin{center}
%%\epsfig{file=Fig1a.eps,width=1\linewidth}\\
%\epsfig{file=Fig1a.eps,width=0.6\linewidth}\\
%%\def\svgwidth{1\linewidth}
%%\input{fig1a.eps_tex}
%%\epsfig{file=fig1a,width=1\linewidth}\\
%\vspace{0.5cm}
%%
%\def\svgwidth{.3\linewidth}
%\input{fig1b.eps_tex}
%%\epsfig{file=fig1b,width=4.5cm}
%\end{center}
%\vspace{-0.3cm}
%\caption{Construction of tensors. {\em Top:} Tensorization of a vector or matrix into the so-called {\em quantized format}; in scientific computing this facilitates super-compression of large-scale vectors or matrices. %(e.g., through a polyadic or a TT decomposition).
%{\em Bottom:} Tensor formed through the discretization of a trivariate function $f(x,y,z)$.}
%\label{Fig:Tensorization}
%\end{figure}

%\begin{figure}[t!]
%%\label{Fig:Tensorization}
%\begin{center}
%%\epsfig{file=Tensvect-CPD.eps,width=1\linewidth}\\
%\def\svgwidth{1\linewidth}
%\input{fig1a.eps_tex}
%%\epsfig{file=fig1a,width=1\linewidth}\\
%\vspace{0.5cm}
%%
%\def\svgwidth{.5\linewidth}
%\input{fig1b.eps_tex}
%%\epsfig{file=fig1b,width=4.5cm}
%\end{center}
%\vspace{-0.3cm}
%\caption{Construction of tensors. {\em Top:} Tensorization of a vector or matrix into the so-called {\em quantized format}; in scientific computing this facilitates super-compression of large-scale vectors or matrices. %(e.g., through a polyadic or a TT decomposition).
%{\em Bottom:} Tensor formed through the discretization of a trivariate function $f(x,y,z)$.}
%\label{Fig:Tensorization}
%\end{figure}

The procedure of creating a data tensor from lower-dimensional original data is referred to as {\em tensorization}, and we propose the following taxonomy for tensor generation:
\begin{itemize}
\item[1)] {\em Rearrangement of lower dimensional data structures}.  Large-scale vectors or matrices are readily tensorized to higher-order tensors, and can be compressed through tensor decompositions if they admit a low-rank tensor approximation; this principle facilitates big data analysis  \cite{Khoromskij-TT,Hackbush2012,Grasedyck-rev} (see Figure~\ref{Fig:Tensorization} (top)).
For instance, a one-way exponential signal $x(k) = a z^{k}$ can be rearranged into a rank-1 Hankel matrix or a Hankel tensor \cite{sid_carat}:
\begin{equation}
\mathbf{H} =
\left(\begin{array}{cccc}
x(0) & x(1) & x(2) & \cdots \\
x(1) & x(2) & x(3) & \cdots \\
x(2) & x(3) & x(4) & \cdots  \\
\vdots & \vdots & \vdots
\end{array}
\right) = a \; \bb \circ \bb,
\label{Eq:Hankelization_Exponential}
\end{equation}
 where $\bb= [1,z,z^2,\cdots]^T$.
%\be
%%%\begin{displaymath}
%\mathbf{H} &=&
%\left(\begin{array}{cccc}
%x(0) & x(1) & x(2) & \cdots \\
%x(1) & x(2) & x(3) & \cdots \\
%x(2) & x(3) & x(4) & \cdots  \\
%\vdots & \vdots & \vdots
%\end{array}
%\right) \nonumber \\
%&=&
%a \left( \begin{array}{c}
%1 \\ z \\ z^2 \\ \vdots
%\end{array}
%\right)
%\left( \begin{array}{cccc}
%1 & z & z^2 & \cdots
%\end{array}
%\right).
%\label{Eq:Hankelization_Exponential}
%%%\end{displaymath}
%\ee

Also, in sensor array processing, tensor structures naturally emerge when combining snapshots from identical subarrays \cite{Sidiropoulos00Bro}.
\item[2)] {\em Mathematical construction}. Among many such examples, the $N$th-order moments (cumulants) of a vector-valued random variable form an $N$th-order tensor \cite{Nikias-Petropulu-HOS}, while in second-order ICA snapshots of data statistics (covariance matrices) are effectively slices of a third-order tensor \cite{SOBI-B97,Comon-Jutten2010}. Also,  a  ({\em channel} $\times$ {\em time}) data matrix can be transformed into a ({\em channel} $\times$ {\em time} $\times$ {\em frequency}) or ({\em channel} $\times$ {\em time} $\times$ {\em scale}) tensor via time-frequency or wavelet representations, a powerful procedure in multichannel EEG analysis in brain science \cite{Miwakeichi,NMF-book}.
\item[3)] {\em Experiment design}. Multi-faceted data can be naturally stacked into a tensor; for instance, in wireless communications the so-called signal diversity (temporal, spatial, spectral, \ldots) corresponds to the order of the tensor  \cite{SidiropoulosGianna}. %\remove{whose modes are the temporal, spatial, spectral etc. components}.
In the same spirit, the standard EigenFaces can be generalized to TensorFaces by combining images with different illuminations, poses, and expressions \cite{Vasilescu}, while the common modes in EEG recordings across subjects, trials, and conditions are best analyzed when combined together into a tensor \cite{Morup11}.
\item[4)] {\em Naturally tensor data}. Some data sources are readily generated as tensors (e.g., RGB color images, videos, 3D light field displays) \cite{tensor-displays}. Also in scientific computing we often need to evaluate a discretized multivariate function; this is  a natural tensor, as illustrated in Figure~\ref{Fig:Tensorization} (bottom) for a trivariate function $f(x,y,z)$ \cite{Khoromskij-TT,Hackbush2012,Grasedyck-rev}. %(http://web.media.mit.edu/$\sim$gordonw/TensorDisplays).
\end{itemize}

The high dimensionality of the tensor format is associated with blessings --- these include possibilities to obtain compact representations, uniqueness of decompositions, flexibility in the choice of constraints, and generality of components that can be identified.

\section{Canonical Polyadic Decomposition} \label{sect:CPD}

\begin{figure*}[t!]
\centering
\includegraphics[width=15.2cm]{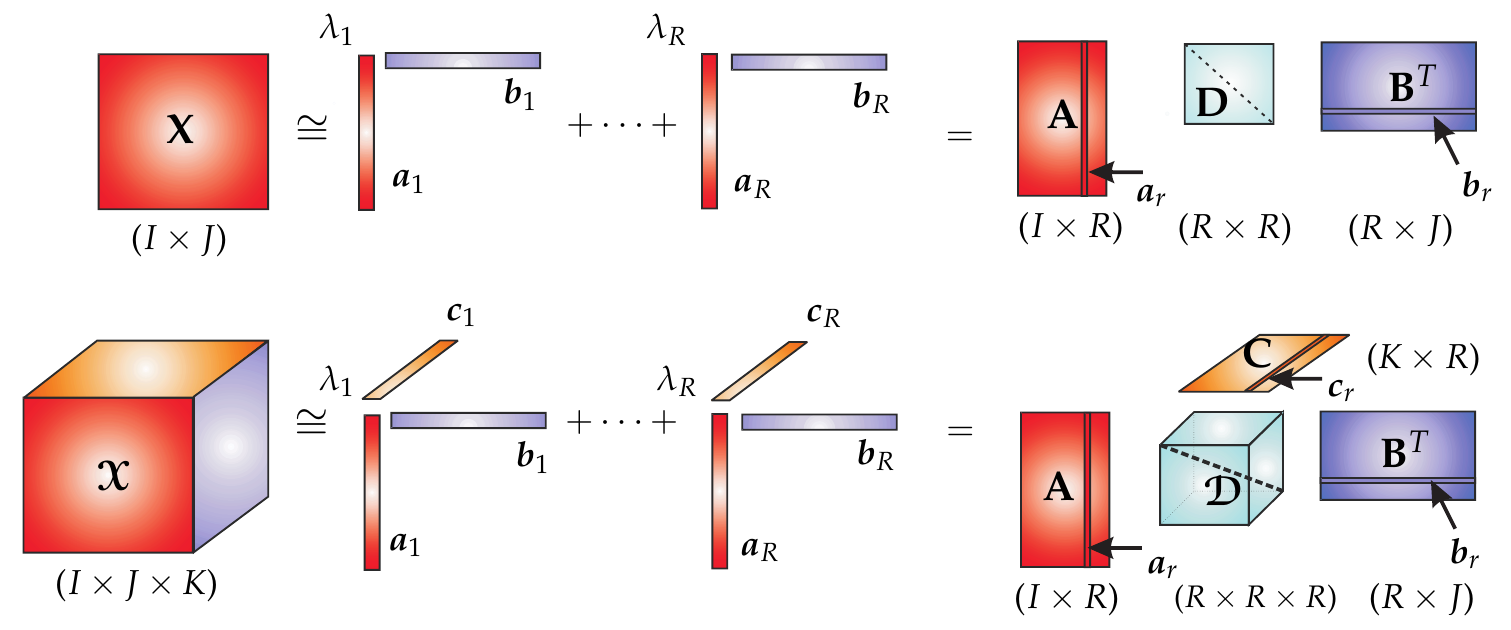}
\caption{Analogy between dyadic (top) and polyadic (bottom) decompositions; the Tucker format has a diagonal core. \replace{Their}{The} uniqueness \replace{properties are}{of these decompositions is}  a prerequisite for blind source separation and latent variable analysis.}
\label{fig:CPD}
\end{figure*}

{\bf Definition.} A Polyadic Decomposition (PD) represents an $N$th-order tensor $\tX \in \Real^{I_1 \times I_2 \times  \cdots \times I_N}$ as a linear combination of  rank-1 tensors in the form
\begin{equation}
\label{CPDN}
 \tX  =  \sum_{r=1}^R  \lambda_r \; \bb^{(1)}_r \circ  \bb^{(2)}_r  \circ \cdots \circ \bb^{(N)}_r.
\end{equation}
Equivalently, $\tX$ is expressed as a multilinear product with a diagonal core:
\be
\label{CPDNpart2}
 \tX & = & \tD \times_1 \bB^{(1)}\times_2 \bB^{(2)} \cdots \times_N \bB^{(N)} \nonumber \\
 & = & \llbracket \tD; \matn[1]{B},  \matn[2]{B}, \ldots, \matn[N]{B}\rrbracket,
\ee
where $\tD = \mbox{\rm diag}_N(\lambda_1, \lambda_2, \ldots, \lambda_R)$ ({\em cf.} the matrix case in (\ref{2CA})).
Figure \ref{fig:CPD} (bottom) illustrates these two interpretations for a third-order tensor. The  tensor rank is  defined as the smallest value of $R$ for which (\ref{CPDN}) holds exactly; the minimum rank PD is called {\emph{canonical}} (CPD) and is desired in  signal separation. The term CPD may also be considered as an abbreviation of CANDECOMP/PARAFAC decomposition, see {\bf Historical notes}.

The matrix/vector form of CPD can be obtained via the Khatri-Rao products as:
\begin{align}
\label{CPD-KR1}
%\bX_{(n)} & = & \bB^{(n)} \bD (\bigodot_{k \neq n} \bB^{(k)})^T \\ \notag \\
\bX_{(n)}  &=  \bB^{(n)} \bD \left(\bB^{(N)} \odot  \cdots \odot \bB^{(n+1)}  \odot \bB^{(n-1)} \odot  \cdots   \odot \bB^{(1)}\right)^T \\
%\bX_{(n)} & = & {\bB^{(n)} \bD (\bB^{(N)} \odot  \cdots \odot \bB^{(n+1)} \odot \bB^{(n-1)} \odot \cdots \odot \bB^{(1)})^T} \notag
\text{vec} (\tX)  &=  [\bB^{(N)} \odot \bB^{(N-1)} \odot \cdots \odot\bB^{(1)}] \,  \bd.\; \; %\notag
\label{CPD-KR2}
\end{align} %}
where  %$\bD = \diag (\lambda_1, \lambda_2, \ldots, \lambda_R)$ and
$\bd = (\lambda_1, \lambda_2, \ldots, \lambda_R)^T$.

%\begin{figure*}[t!]
%\centering
%\def\svgwidth{0.8\linewidth}
%\input{fig2.eps_tex}
%%\psfrag{tX}[t][t]{\scalebox{1}{\color[rgb]{0,0,0}\setlength{\tabcolsep}{0pt}\begin{tabular}{c}$\tX$\end{tabular}}}%
%%\psfrag{tD}[t][t]{\scalebox{1}{\color[rgb]{0,0,0}\setlength{\tabcolsep}{0pt}\begin{tabular}{c}$\tD$\end{tabular}}}%
%%\includegraphics[width=1\textwidth]{fig2_r.eps}
%%\includegraphics[width=7.2cm,height=6.9cm]{SPM-CPTv2.eps}
%%\includegraphics[width=7.2cm,height=6.9cm]{SPM-CPT.eps} \\
%%\includegraphics[width=7.2cm,height=6.9cm]{Fig2b.eps}
%\caption{Analogy between dyadic (top) and polyadic (bottom) decompositions; the Tucker format has a diagonal core. \replace{Their}{The} uniqueness \replace{properties are}{of these decompositions is}  a prerequisite for blind source separation and latent variable analysis.}
%\label{fig:CPD}
%\end{figure*}

{\bf Rank.} As mentioned earlier, rank-related properties are very different for matrices and tensors. For instance, the number of complex-valued rank-1 terms needed to represent a higher-order tensor can be strictly less than the number of real-valued rank-1 terms \cite{Landsberg}, while the determination of tensor rank is in general NP-hard \cite{hastad}. Fortunately, in signal processing applications, rank estimation most often corresponds to determining the number of tensor components that can be retrieved with sufficient accuracy, and often there are only a few data components present. A pragmatic first assessment of the number of components may be through the inspection of the multilinear singular value spectrum (see Section {\bf Tucker Decomposition}), which indicates the size of the core tensor in Figure \ref{fig:CPD} (bottom-right).
Existing techniques for rank estimation include the CORCONDIA algorithm (core consistency diagnostic) which checks whether the core tensor is (approximately) diagonalizable \cite{Smilde},  while
a number of techniques  operate by balancing the approximation error versus the number of degrees of freedom for a varying number of rank-1 terms \cite{Timmerman,Ceulemans,IMM2009-05806}.

{\bf Uniqueness.} Uniqueness conditions give theoretical bounds for exact tensor decompositions.
A classical uniqueness  condition  is due to Kruskal \cite{Kruskal1977}, which states that for third-order tensors
the CPD is unique up to unavoidable scaling and permutation ambiguities, provided that $k_{\bB^{(1)}} + k_{\bB^{(2)}} + k_{\bB^{(3)}} \geq 2R + 2$, where the Kruskal rank $k_{\bB}$ of a matrix $\bB$ is the maximum value ensuring that any subset of $k_{\bB}$ columns is linearly independent.
In sparse modeling, the  term $(k_B +1 )$ is also known as \add{the} spark \cite{bruckstein2009sparse}.
A generalization to $N$th-order tensors is due to Sidiropoulos and Bro \cite{Sidiropoulos2000}\add{, and is given by:}
\be
\sum_{n=1}^N k_{\bB^{(n)}} \geq 2R + N-1.
\ee
More relaxed uniqueness conditions can be obtained when one factor matrix has full column rank  \cite{JiangSidiropoulos,Lathauwer06-rank,Stegman09}; for a thorough study \replace{in}{of} the third-order case, we refer to \cite{Domanov12}.
This \add{all} shows that, compared to matrix decompositions, CPD is unique under more natural and relaxed conditions, that only require the components to be  ``sufficiently different'' and  their number not unreasonably large. These conditions do not have a matrix counterpart, and are at the heart of tensor based signal separation.

{\bf Computation.}
Certain conditions, including Kruskal's, enable explicit computation of the factor matrices in (\ref{CPDN}) using linear algebra (essentially, by solving sets of linear equations and by \replace{computation of}{computing} (generalized) Eigenvalue Decomposition) \cite{Harshman,Sanchez1990,Lathauwer06-rank,DomanovredGEVD}. The presence of noise in data means that CPD is rarely exact, and we need to fit a CPD model to the data by minimizing a suitable cost function. %The noise-free solution may be used for initialization at high SNR.
This is typically achieved by minimizing the Frobenius norm of the difference between the given data tensor and its CP approximation, or alternatively by least absolute error fitting when the noise is Laplacian \cite{Vorobyov}. Theoretical Cram\'er-Rao Lower Bound (CRLB) and Cram\'er-Rao Induced Bound (CRIB) for \add{the assessment of} CPD performance were %first
derived %by Liu and Sidiropoulos
in \cite{Liu01cramer-raolower} and \cite{Petr_CRIB}.

Since the computation of CPD is intrinsically multilinear, we can arrive at the solution through a sequence of linear sub-problems as in the Alternating Least Squares (ALS) framework, whereby  the LS cost function is optimized for one component matrix at a time, while keeping the other component matrices fixed \cite{Harshman}. As seen  from (\ref{CPD-KR1}), such a conditional update scheme boils down to solving overdetermined sets of linear equations. %\cite{phanref42}.

While the ALS is attractive for its simplicity and satisfactory performance for a few well separated components and at sufficiently high SNR, it also inherits the problems of alternating algorithms and is not guaranteed to converge to a stationary point.
This can be rectified by only updating the factor matrix for which the cost function has most decreased at a given step \cite{chen2012}, but this results in an $N$-times increase in computational cost per iteration. The convergence of ALS is not yet completely understood --- it is quasi-linear close to the stationary point \cite{uschmajew2012}, while it becomes rather slow for ill-conditioned cases; for more detail we refer to \cite{mohlenkamp2013,razaviyayn2013unified}.

Conventional all-at-once algorithms for numerical optimization such as nonlinear conjugate gradients, quasi-Newton or nonlinear least squares \cite{Paatero99,Acar-Morup11} have been shown to often outperform ALS for ill-conditioned cases and to be typically more robust to overfactoring, but come at a cost of a much higher computational load per iteration. More sophisticated versions use the rank-1 structure of the terms within CPD to perform efficient computation and storage
of the Jacobian and (approximate) Hessian; their complexity is on par with ALS while for ill-conditioned cases the performance is often superior  \cite{Phan2012-Hess,Sorber2012b}.

An important difference between matrices and tensors is that the existence of a best rank-$R$ approximation of a tensor of rank greater than $R$ is not guaranteed { \cite{Landsberg,deSilva-Lim08}} since the set of tensors whose rank is at most $R$ is not closed. As a result, cost functions for computing factor matrices may only have an infimum (instead of a minimum) so that their minimization will approach the boundary of that set without ever reaching the boundary point.
This will cause two or more rank-1 terms go to infinity upon convergence of an algorithm, however, numerically the diverging terms will almost completely cancel one another  while  the overall cost function will still decrease along the iterations \cite{krijnen}.
These diverging terms indicate an inappropriate data model: the mismatch between the CPD and the original data tensor may arise due to an underestimated number of components, not all tensor components having a rank-1 structure, or  data being too noisy.

{\bf Constraints.} As mentioned earlier, under quite mild conditions the CPD is unique by itself, without requiring additional constraints. However, in order to enhance \add{the} accuracy and robustness with respect to noise, prior knowledge of data properties (e.g., statistical independence, sparsity) may be incorporated into \add{the} constraints on factors  so as to facilitate \add{their} physical interpretation\remove{of the factors}, relax the uniqueness conditions, and even simplify computation \cite{sorensen,sorensenVDM,Zhou2012-SPL}. Moreover, the orthogonality and nonnegativity constraints ensure the existence of the minimum of the optimization criterion \add{used} \cite{krijnen,limcomNN,sorensen}.

{\bf Applications.} The CPD has already been established as an advanced tool for signal separation in vastly diverse branches of signal processing and data analysis, such as in audio and speech processing, biomedical engineering, chemometrics, and machine learning \cite{Smilde,Acar2008,Kolda08,Morup11}.
Note that algebraic ICA algorithms are effectively based on the CPD of a tensor of the statistics of recordings; the statistical independence of the sources is reflected in the diagonality of the core tensor in Figure  \ref{fig:CPD}, that is, in vanishing cross-statistics \cite{comon94,Comon-Jutten2010}.
The CPD is also heavily used in exploratory data analysis, where the rank-1 terms capture %\remove{the}
essential properties of dynamically complex signals \cite{Kroonenberg}.
Another example is in wireless communication, where the signals transmitted by different users correspond to rank-1 terms in the case of line-of-sight propagation \cite{Sidiropoulos00Bro}.
%; then the data streams are conveniently stacked into a tensor if the system has at least triple diversity  \cite{Sidiropoulos00Bro}.
Also, in harmonic retrieval and direction of arrival type applications, real or complex exponentials have a rank-1 structure, for which the use of CPD is natural \cite{sid_carat,sorensenVDM}.

{\bf Example 1.}
Consider a sensor array consisting of $K$ displaced but otherwise identical subarrays of $I$ sensors, with $\tilde{I} = KI$ sensors in total. For $R$ narrowband sources in the far field, the baseband equivalent model of the array output becomes $\bX = \bA \bS^T + \bE$, where $\bA \in \mathbb{C}^{\tilde{I} \times R}$ is the global array response, $\bS \in \mathbb{C}^{J \times R}$ contains $J$ snapshots of the sources, and $\bE$ is noise.
A single source ($R=1$) can be obtained from the best rank-1 approximation of the matrix $\bX$, however, for $R> 1$ the decomposition of $\bX$ is not unique, and hence the separation of sources is not possible without incorporating additional information. Constraints on the sources that may yield a unique solution are, for instance, constant modulus or statistical independence \cite{Veen,Comon-Jutten2010}.

Consider a row-selection matrix $\bJ_k \in \mathbb{C}^{I \times \tilde{I}}$ that extracts the rows of $\bX$ corresponding to the $k$-th subarray, $k=1, \ldots, K$. %, and consider $\bJ_k \bX \in \mathbb{C}^{I \times J}$.
For two identical subarrays, the generalized EVD  of the matrices $\bJ_1 \bX$ and $\bJ_2 \bX$ corresponds to the well-known ESPRIT \cite{roy1989esprit}. For the case $K > 2$, we shall consider $\bJ_k \bX$ as slices of the tensor $\tX \in \mathbb{C}^{I \times J \times K}$ (see Section {\bf Tensorization}). It can be shown that the signal part of $\tX$ admits a CPD as in (\ref{CPDN})--(\ref{CPDNpart2}), with $\lambda_1 = \cdots = \lambda_R = 1$, $\bJ_k \bA = \bB^{(1)} \mbox{\rm diag}(\bb^{(3)}_{k1}, \ldots, \bb^{(3)}_{kR})$ and $\bB^{(2)} = \bS$ \cite{Sidiropoulos00Bro}, and the consequent source separation under rather mild conditions --- its uniqueness does not require constraints such as statistical independence or constant modulus. Moreover, the decomposition is unique even in cases when the number of sources $R$ exceeds the number of subarray sensors $I$, or even the total number of sensors $\tilde{I}$.
Notice that particular array geometries, such as linearly and uniformly displaced subarrays, can be converted into a constraint on CPD, yielding \add{a} further \replace{relaxed}{relaxation of the} uniqueness conditions, reduced sensitivity to noise, and often faster computation \cite{sorensenVDM}.

\section{Tucker Decomposition} \label{sect:Tucker}

\begin{figure}[t]
    \begin{center}
\includegraphics[width=9.9cm]{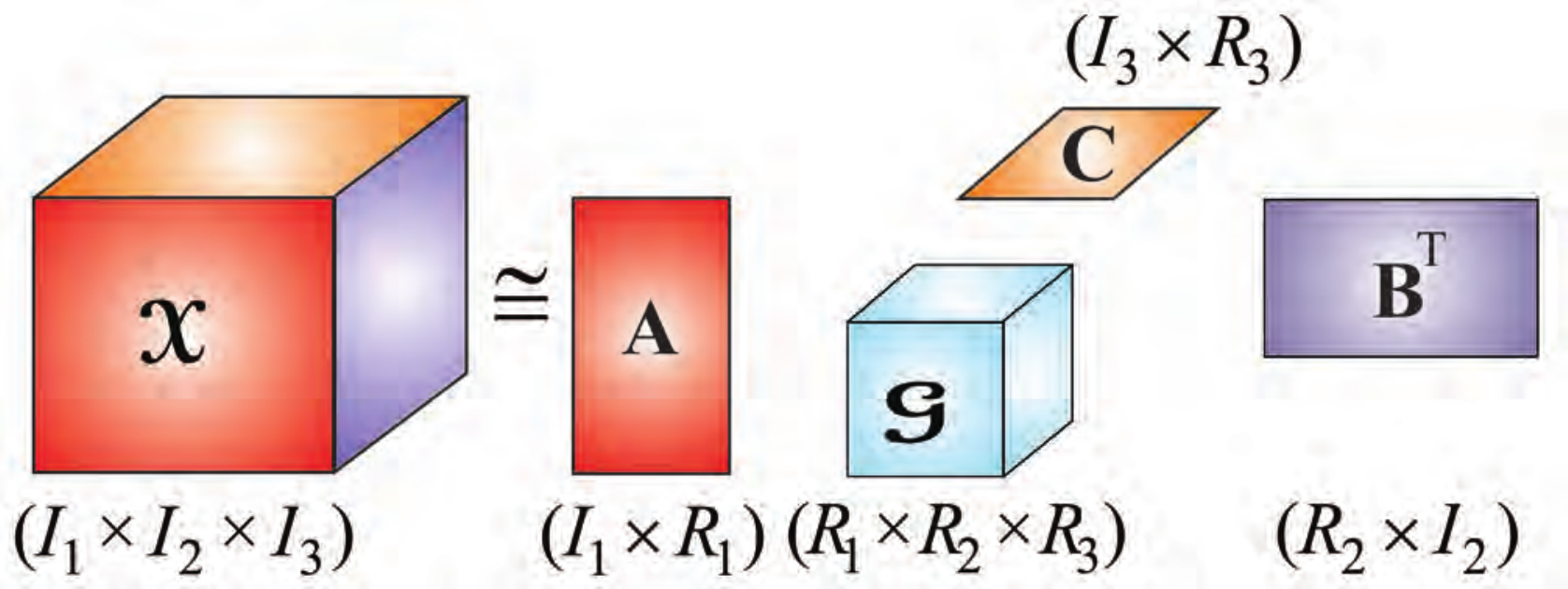}
    \end{center}
\caption{Tucker decomposition of a third-order tensor. The column spaces of $\bA$, $\bB$, $\bC$ represent the signal subspaces for the three modes. The core tensor $\tG$ is nondiagonal, accounting for possibly complex interactions among tensor components. }
\label{fig:Tucker3bis}
\end{figure}

Figure \ref{fig:Tucker3bis} illustrates the principle of Tucker decomposition which treats a tensor $\tX \in \Real^{I_1\times I_2 \times \cdots \times I_N}$ as a multilinear transformation of a (typically dense but small) core tensor $\tG  \in \Real^{R_1 \times R_2 \times \cdots \times R_N}$ by the factor matrices $\bB^{(n)} =[\bb_1^{(n)},\bb_2^{(n)},\ldots,\bb_{R_n}^{(n)}] \in \Real^{I_n\times R_n}$, $n=1,2,\ldots,N$\cite{tucker64extension,Tucker1966}, given by
\begin{equation}
\label{eq:Tuckersumofr1}
  {\tX}
  = \sum\limits_{r_1 = 1}^{R_1} \sum\limits_{r_2 = 1}^{R_2} {
 {\cdots
 \sum\limits_{r_N = 1}^{R_N}{g_{r_1 r_2 \cdots r_N} \left(\bb^{(1)}_{r_1} \circ \bb^{(2)}_{r_2} \circ \cdots \circ \bb^{(N)}_{r_N} \right)}}}
 %\notag \\
\end{equation}
or equivalently
\be
 \label{GeneralTDModel}
  \tX & = & \tG \times_1 \bB^{(1)}\times_2 \bB^{(2)} \cdots \times_N \bB^{(N)} \nonumber \\
  &=& \llbracket\tG; \bB^{(1)},\bB^{(2)},\ldots,\bB^{(N)}\rrbracket.
 \ee
Via the Kronecker products (see Table~\ref{table_notation2}) Tucker decomposition can be expressed in a  matrix/vector form as:
\be
%\begin{align}
\label{Tucker-kron}
%\bX_{(n)} & =  \bB^{(n)} \bG_{(n)} (\bigotimes_{k \neq n} \bB^{(k)})^T
 \bX_{(n)} & =&  \bB^{(n)} \bG_{(n)} (\bB^{(N)} \otimes \cdots \otimes \bB^{(n+1)}
  \otimes \; \bB^{(n-1)} \otimes \cdots \otimes \bB^{(1)})^T  \\ %\notag \\
 \text{vec}(\tX) & = & [\bB^{(N)} \otimes \bB^{(N-1)} \otimes \cdots \otimes \bB^{(1)}] \; \text{vec}(\tG). \qquad
%\end{align}
\ee
%%\be
%\begin{align}
%\label{Tucker-kron}
%%\bX_{(n)} & =  \bB^{(n)} \bG_{(n)} (\bigotimes_{k \neq n} \bB^{(k)})^T
%& \bX_{(n)}  =  \bB^{(n)} \bG_{(n)} (\bB^{(N)} \otimes \cdots \otimes \bB^{(n+1)}  \otimes \bB^{(n-1)} \otimes \cdots \otimes \bB^{(1)})^T \notag \\
%& \text{vec}(\tX)  =  [\bB^{(N)} \otimes \bB^{(N-1)} \otimes \cdots \otimes \bB^{(1)}] \; \text{vec}(\tG). \qquad \notag
%\end{align}
%
Although Tucker initially used the orthogonality and ordering constraints on the core tensor and factor matrices \cite{tucker64extension,Tucker1966}, we can also employ other meaningful constraints (see below).
 %\cite{tenberge-sp,Kiers1998e,Rocci:2002wr}.

% \begin{figure}[t]
%    \begin{center}
%    \def\svgwidth{1\linewidth}
%\input{Fig3.eps_tex}
%%\includegraphics[width=7.2cm]{fig3}
%    \end{center}
%\caption{Tucker decomposition of a third-order tensor. The column spaces of $\bA$, $\bB$, $\bC$ represent the signal subspaces for the three modes. The core tensor $\tG$ is nondiagonal, accounting for possibly complex interactions among tensor components. }
%\label{fig:Tucker3bis}
%\end{figure}

{\bf Multilinear rank.} For a core tensor of minimal size, $R_1$ is  the column rank (the dimension of the subspace spanned by mode-1 fibers), $R_2$ is the row rank (the dimension of the subspace spanned by mode-2 fibers), and so on. A remarkable difference from matrices is that the values of $R_1,R_2,\ldots,R_N$ can be  different for $N\geq3$. The $N$-tuple $(R_1,R_2,\ldots,R_N)$ is consequently called the {\emph{multilinear rank}} of the tensor $\tX$.

{\bf Links between CPD and Tucker decomposition.}
Eq. (\ref{eq:Tuckersumofr1}) shows that Tucker decomposition can be considered as an expansion in rank-1 terms (polyadic but not necessary canonical), while (\ref{CPDNpart2}) represents CPD as a multilinear product of a core tensor and factor matrices (but the core is not necessary minimal); Table \ref{CPD-Tucker} shows various other connections. However, despite the obvious interchangeability of notation, the CP and Tucker decompositions serve different purposes. In general, the Tucker core cannot be diagonalized, while the number of CPD terms may not be bounded by the multilinear rank. Consequently, in signal processing and data analysis, CPD is typically  used for factorizing data into easy to interpret components (i.e., the rank-1 terms), while the goal of unconstrained Tucker decompositions is most often to compress data into a tensor of smaller size (i.e., the core tensor) or to find the subspaces spanned by the fibers (i.e., the column spaces of the factor matrices).

\begin{table}[h!]
\caption{Different forms of CPD and Tucker representations of a third-order tensor $\tX \in \Real^{I \times J \times K}$.}
 \centering
 {\footnotesize  \shadingbox{
    \begin{tabular*}{.9\textwidth}[t]{@{\extracolsep{\fill}}l@{\hspace{1em}}l} \hline  &  \\
{ \raisebox{0mm}[0mm][4mm]{CPD}}
&  Tucker Decomposition \\ \hline
 &   \\
 \multicolumn{2}{c}{ Tensor representation, outer products}\\
 & \\
%\begin{minipage}[t]{.5\textwidth}
\begin{tabular*}{.6\linewidth}[t]{@{\extracolsep{\fill}}cc@{}}
 $\;\;\tX = \sum\limits_{r=1}^{R} {\lambda_r \; \ba_{r} \circ \bb_{r} \circ \bc_{r}} $  & \hspace{6.8cm}
$\tX = \sum\limits_{r_1=1}^{R_1} \sum\limits_{r_2=1}^{R_2} \sum\limits_{r_3=1}^{R_3} {g_{r_1 \,r_2 \, r_3} \; \ba_{r_1} \circ \bb_{r_2} \circ \bc_{r_3}} $
\end{tabular*}
%\end{minipage}
\\
 &  \\ \hline
  & \\
 \multicolumn{2}{c}{Tensor representation, multilinear products}\\
  & \\
   $ \;\; \tX = \tD \times_1 \bA \times_2 \bB \times_3 \bC $  &  \hspace{0em}
 $ \tX = \tG \times_1 \bA \times_2 \bB \times_3 \bC $ \\  &   \\ \hline
    % & \\
% {} & {\hspace{-1.4cm} Tensor representation, multilinear products}\\
%  & \\
%   $ \tX = {\rm diag}_3(\lambda_1, \lambda_2, \ldots, \lambda_R) \times_1 \bA \times_2 \bB \times_3 \bC $  &
% $ \tX = \tG \times_1 \bA \times_2 \bB \times_3 \bC $ \\  &   \\ \hline
     & \\
 \multicolumn{2}{c}{ Matrix representations}\\
    &  \\
       $\;\;\bX_{(1)} = \bA \; \bD \; (\bC \odot \bB)^T $  &
       $\bX_{(1)} = \bA \;\bG_{(1)} \; (\bC \otimes \bB)^T $ \\
      {$\;\;\bX_{(2)} = \bB \; \bD \; (\bC \odot \bA)^T $} &
      {$\bX_{(2)} = \bB \;\bG_{(2)} \; (\bC \otimes  \bA)^T $} \\
       {$\;\;\bX_{(3)} = \bC \; \bD \; (\bB\odot \bA)^T $}
        & {$\bX_{(3)} =\bC \;\bG_{(3)} \; (\bB \otimes  \bA)^T $ } \\
       &  \hspace{-6em} \\ \hline
    & \\
    \multicolumn{2}{c}{Vector representation} \\
                    & \\
                       $\;\;\text{vec}(\tX) =  (\bC \odot \bB \odot \bA) \mbi \bd$
                    % \mbi \lambda =[\lambda_1, \lambda_2,\ldots,\lambda_R]^T$$
                      &  \hspace{-2em}
                       {$\text{vec}(\tX) = (\bC \otimes  \bB \otimes  \bA) \;
                        \text{vec}(\tG) $} \\  &  \hspace{-6em}   \\ \hline
                      & \\
  \multicolumn{2}{c}{ Scalar representation}\\
 & \\
 \begin{tabular*}{.5\linewidth}[t]{@{\extracolsep{\fill}}cc@{}}
     $ \;\;x_{ijk} = \sum\limits_{r=1}^{R} {\lambda_r \,a_{i\,r} \, b_{j\,r} \, c_{k\,r}}  $  & \hspace{6.8cm} %\hspace{-4em}
$ x_{ijk} =\sum\limits_{r_1=1}^{R_1} \sum\limits_{r_2=1}^{R_2} \sum\limits_{r_3=1}^{R_3} {g_{r_1 \, r_2 \, r_3} \,a_{i\,r_1} \, b_{j\,r_2} \, c_{k\,r_3}}  $
\end{tabular*}
                   \\  &   \\ \hline
%    &  \\
   % $\text{vec}(\tX)  \cong (\bC \odot \bB \odot \bA) \;{\boldsymbol \lambda} $ & {$\text{vec}(\tX) \cong (\bC \otimes  \bB \otimes  \bA) \;\text{vec}(\tG) $} \\  &   \\ \hline
       & \\
  \multicolumn{2}{c}{Matrix slices $\bX_k = \bX(:,:,k)$} \\
    &  \\
    $\;\;\bX_k =  \bA \, {\rm diag}(c_{k1}, c_{k2}, \ldots, c_{kR}) \, \bB^T$  &
      $\bX_k= \bA\, \sum\limits_{r_3=1}^{R_3} c_{k r_3} \bG(:,:,r_3) \, \bB^T $ \\
      &  \\ \hline
    \end{tabular*}
    }}
\label{CPD-Tucker}
\end{table}

{\bf Uniqueness.}  The unconstrained Tucker decomposition is in general
not unique,  that is, factor matrices $\bB^{(n)}$ are rotation invariant. However, physically, the subspaces defined by the factor matrices in Tucker decomposition are unique, while the bases in these subspaces may be chosen arbitrarily  --- their choice is compensated for within the core tensor.
This becomes clear upon realizing that any factor matrix in  (\ref{GeneralTDModel})  can be post-multiplied by any  nonsingular (rotation) matrix; in turn, this multiplies the  core tensor by its inverse, that is
\be
  {\tX}
&=& \llbracket \tG; \bB^{(1)},\bB^{(2)}, \ldots, \bB^{(N)} \rrbracket \notag \\
&=&
\llbracket \tH;  \bB^{(1)} \bR^{(1)}, \bB^{(2)}  \bR^{(2)},\ldots, \bB^{(N)} \bR^{(N)} \rrbracket, \notag \\
\tH &=& \llbracket \tG ; \bR^{(1)^{-1}} , \bR^{(2)^{-1}} , \ldots, \bR^{(N)^{-1}} \rrbracket,
\ee
where $\bR^{(n)}$ are invertible.

\begin{figure*}[ht!]
%\label{Fig:MWCA}
\begin{center}
\includegraphics[width=14.5cm,height=5.5cm]{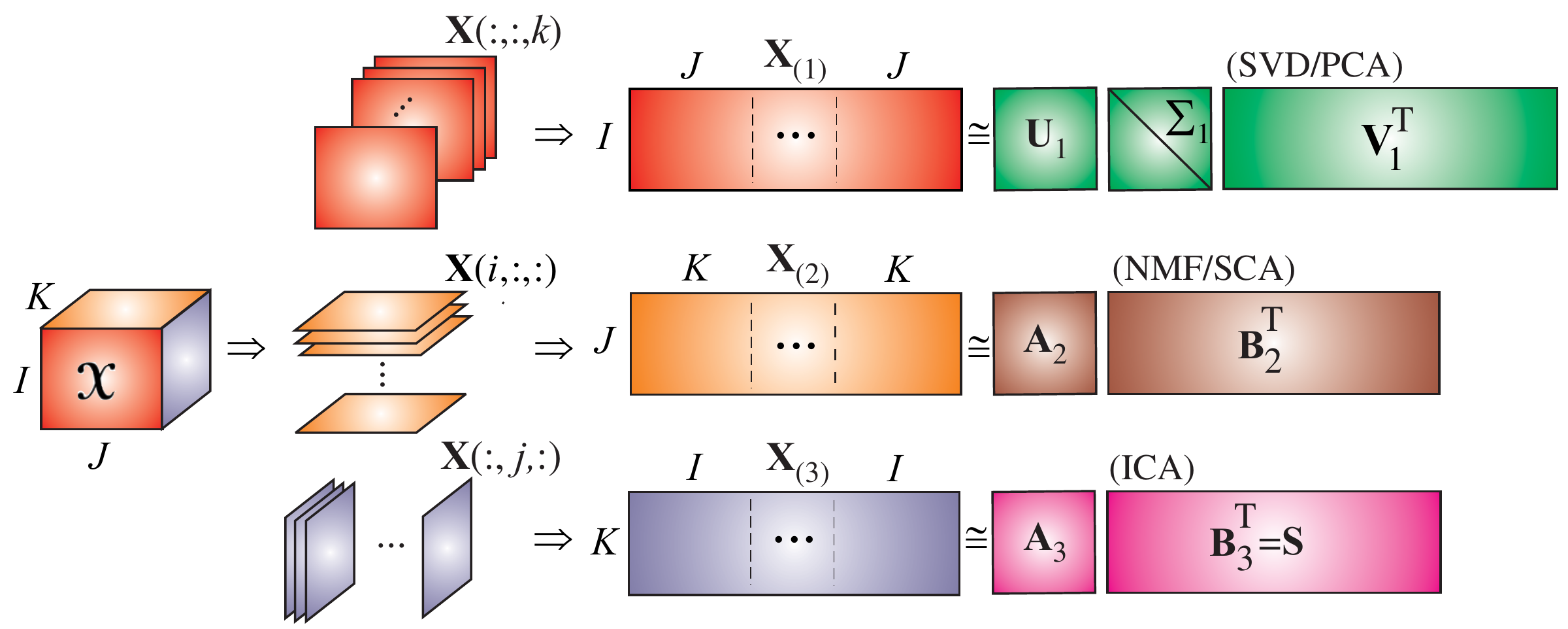}
\end{center}
\caption{Multiway Component Analysis (MWCA) for a third-order tensor, assuming that \add{the} components are: principal and orthogonal
 in the first mode, nonnegative and sparse in the second mode and  statistically independent in the third mode.} %This flexibility makes MWCA more versatile than multilinear SVD and multilinear ICA \cite{MICA2005}. }
\label{Fig:MWCA}
\end{figure*}
%
%\begin{figure*}[ht!]
%%\label{Fig:MWCA}
%\begin{center}
%\def\svgwidth{.9\linewidth}
%\input{Fig4.eps_tex}
%%\epsfig{file=MWCA.eps,width=11cm,height=6}
%%\includegraphics[width=\linewidth,height=4.5cm]{fig4_r}
%\end{center}
%\caption{Multiway Component Analysis (MWCA) for a third-order tensor, assuming that \add{the} components are: principal and orthogonal
% in the first mode, nonnegative and sparse in the second mode and  statistically independent in the third mode.} %This flexibility makes MWCA more versatile than multilinear SVD and multilinear ICA \cite{MICA2005}. }
%\label{Fig:MWCA}
%\end{figure*}

{\bf Multilinear SVD (MLSVD).} %
Orthonormal bases in a constrained  Tucker representation can be obtained via the SVD of the mode-$n$ matricized tensor $\bX_{(n)} = \bU_{n} \mbi {\Sigma}_{n} \bV_{n}^T$ (i.e., $\bB^{(n)}=\bU_{n}$, $n=1,2,\ldots,N$).
Due to the orthonormality, the corresponding core tensor becomes
\be
\tS = \tX  \times_1 \bU_{1}^T  \times_2 \bU_{2}^T  \cdots  \times_N \bU_{N}^T. %,
\label{HOSVD:core}
\ee

Then, the singular values of $\bX_{(n)}$ are the Frobenius norms of the corresponding slices of the core tensor $\tS$: $({\Sigma}_{n})_{r_n,r_n} = \| \tS(:,:,\ldots,r_n,:,\ldots,:) \|$, with slices in the same mode being mutually orthogonal, %w.r.t. the tensor inner product,
i.e., their inner products are zero.
The columns of $\bU_{n}$ may thus be seen as multilinear singular vectors, while the norms of the slices of the core are multilinear singular values \cite{Lathauwer01}. As in the matrix case, the multilinear singular values govern the multilinear rank, while the multilinear singular vectors allow, for each mode separately, an interpretation as in PCA \cite{Kroonenberg}.

{\bf Low multilinear rank approximation.}  Analogous to PCA, a large-scale data tensor $\tX$ can be approximated by discarding the multilinear singular vectors and slices of the core tensor that correspond to small multilinear singular values, that is, through truncated  matrix SVDs.
Low multilinear rank approximation is always well-posed, however, the truncation is not necessarily optimal in the LS sense, although a good estimate can often be made as the approximation error corresponds to the degree of truncation. When it comes to finding the best approximation, the ALS type algorithms exhibit similar advantages and drawbacks to those  used for CPD \cite{Lathauwer_HOOI,Kroonenberg}. Optimization-based algorithms exploiting second-order information have also been proposed \cite{Savas2010,IshtevaAHL11}.

{\bf Constraints and Tucker-based multiway component analysis (MWCA).}
Besides orthogonality, constraints that may help to find unique basis vectors in \add{a} Tucker representation include statistical independence, sparsity, smoothness and nonnegativity \cite{NMF-book,Zhou-Cichocki-MBSS,Cichocki-GCA}. Components of a data tensor  seldom have the same properties in its modes, and for physically meaningful representation different constraints may be required in different modes, so as to match the properties of the data at hand. Figure~\ref{Fig:MWCA}
illustrates the concept of MWCA and its flexibility in choosing the mode-wise constraints; \remove{while} a Tucker representation of MWCA naturally accommodates such diversities in different modes.

{\bf Other applications.}
We have shown that Tucker decomposition may be considered as a multilinear extension of PCA \cite{Kroonenberg}; it therefore generalizes signal subspace techniques, with applications including classification, feature extraction, and subspace-based harmonic retrieval \cite{Haardt08,Vasilescu,Phan2010TF,Lu-2011}. For instance, a low multilinear rank approximation achieved through
Tucker decomposition may yield a higher Signal-to-Noise Ratio (SNR) than the SNR in the original raw data tensor, making Tucker decomposition a very natural tool for compression and signal enhancement \cite{Kroonenberg,Smilde,Comon-ALS09}.

\section{Block Term Decompositions}

\begin{figure}[ht]
\centering
\includegraphics[width=13.2cm]{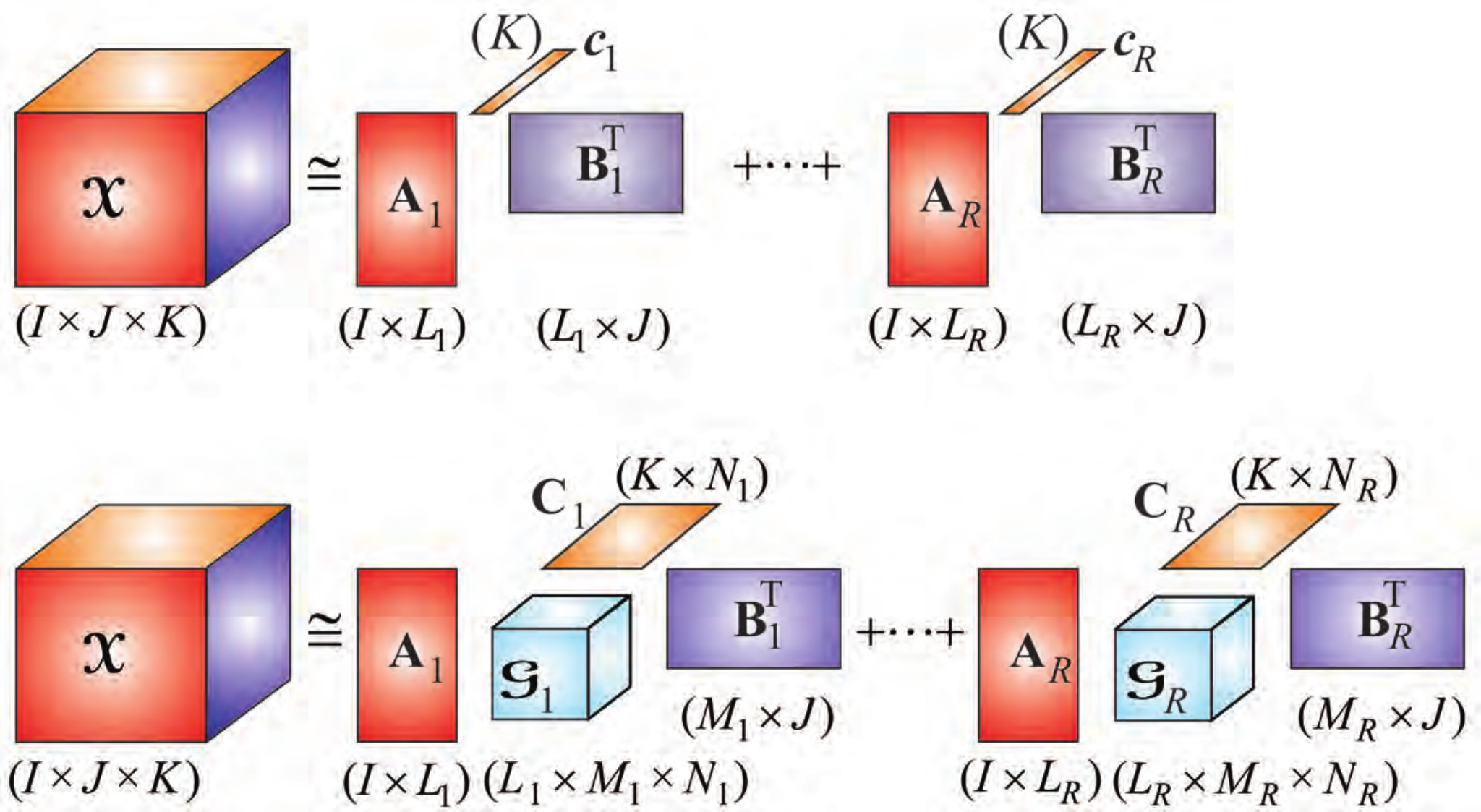}
\caption{Block Term Decompositions (BTDs) find data components that are structurally more complex than the rank-1 terms in CPD.
%$\tX \in \Real^{I \times J \times K}$:
{\em Top:} Decomposition into terms with multilinear rank $(L_r,L_r,1)$. {\em Bottom:} Decomposition into terms with multilinear rank $(L_r,M_r,N_r)$.}
\label{fig-BTD}
\end{figure}

We have already shown that CPD is unique under quite mild conditions, a further advantage of tensors over matrices is that it is even possible to relax the rank-1 constraint on the terms, thus opening completely new possibilities in e.g. BSS. For clarity, we shall consider the third-order case, whereby, by replacing the rank-1 matrices $\bb^{(1)}_r \circ  \bb^{(2)}_r = \bb^{(1)}_r   {\bb^{(2)\; T}_r}$ in (\ref{CPDN}) by low-rank matrices $\bA_r \bB_r^T$, the tensor $\tX$ can be represented as (Figure \ref{fig-BTD}, top):
\begin{equation}
\tX =  \sum_{r=1}^R (\bA_r \bB_r^T) \circ \bc_r.
\label{eq:BTD-LL1}
\end{equation}
Figure \ref{fig-BTD} (bottom) shows that we can even use terms that are only required to have a low multilinear rank (see also Section {\bf Tucker Decomposition}), to give:
\begin{equation}
\tX =  \sum_{r=1}^R \tG_r  \times_1 \bA_r  \times_2 \bB_r  \times_3 \bC_r.
\label{eq:BTD-LMN}
\end{equation}
These so-called Block Term Decompositions (BTD) admit the modelling of more complex signal components than CPD, and are unique under more restrictive but still fairly natural conditions \cite{Lath-BCM12,LathauwerTBSS,LathauwerLVA2012}.

%\begin{figure}[ht]
%\centering
%\def\svgwidth{1\linewidth}
%\input{Fig5a.eps_tex}
%%\includegraphics[width=7.2cm,height=2.5cm]{fig5a} \\
%\vspace{0.3cm}
%%\hspace{-0.2cm}\includegraphics[width=7.4cm,height=2.7cm]{fig5b_r1}
%\def\svgwidth{1\linewidth}
%\input{Fig5b.eps_tex}
%\vspace{.1ex}
%\caption{Block Term Decompositions (BTDs) find data components that are structurally more complex than the rank-1 terms in CPD.
%%$\tX \in \Real^{I \times J \times K}$:
%{\em Top:} Decomposition into terms with multilinear rank $(L_r,L_r,1)$. {\em Bottom:} Decomposition into terms with multilinear rank $(L_r,M_r,N_r)$.}
%   %(for simplicity, we assumed that all core tensors
%   %$\tG_r \in \Real^{R_1 \times R_2 \times R_3}$ have the same dimensions).}
%\label{fig-BTD}
%\end{figure}
%%

\begin{figure}[ht]
   % \centering
     \includegraphics[width=16.9cm]{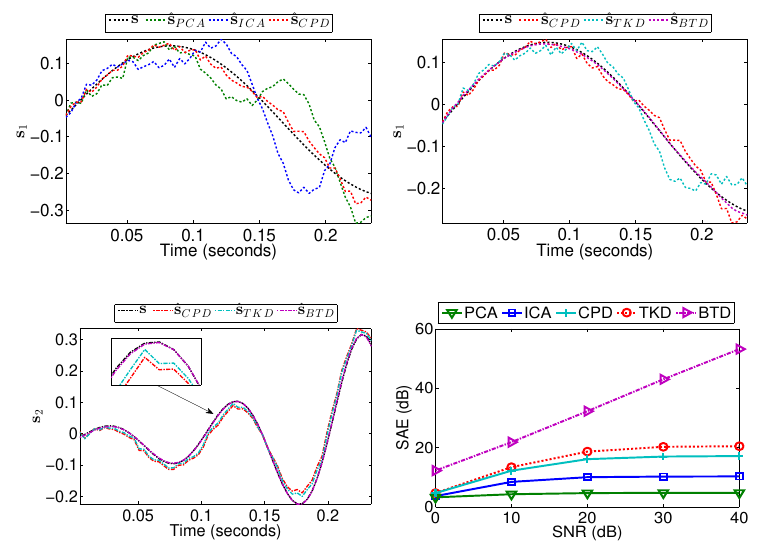}
     \captionsetup{singlelinecheck=off}
\caption[foo bar]{Blind separation of the mixture of a pure sine wave and an exponentially modulated sine wave using PCA, ICA, CPD, Tucker decomposition (TKD)
and BTD. The sources $s_1$ and $s_2$ are correlated and of  short duration; the symbols $\hat s_1$ and $\hat s_2$ denote the estimated sources.}
\label{fig_ex_sine_exp}
\end{figure}

{\bf Example 3.}
To compare some standard and tensor approaches for the separation of short duration correlated sources, BSS was performed on five linear mixtures of the sources
$s_1(t) = \sin(6\pi t)$ and
$s_2(t) = \exp(10t)   \sin(20 \pi t )$,
which were contaminated by white Gaussian noise, to give the mixtures $\bX = \bA \bS + \bE \in \Real^{5 \times 60}$, where $\bS(t) = [s_1(t), s_2(t)]^T$ and $\bA \in \Real^{5 \times 2}$ was a random matrix whose columns (mixing vectors) satisfy $\ba_1^T \ba_2 = 0.1$, $\|\ba_1\|= \|\ba_2\| = 1$.
The 3Hz sine wave did not complete a full period over the 60 samples, so that the two sources had a correlation degree of $\frac{|\bs_1^T \bs_2|}{\|\bs_1\|_2 \|\bs_2\|_2} =  0.35$. The tensor approaches, CPD, Tucker decomposition and BTD employed a third-order tensor $\tX$ of size $24 \times  37  \times  5$ generated from five Hankel matrices whose elements obey $\tX(i,j,k) = \bX(k,i+j-1)$ (see Section {\bf Tensorization}). %The plots compare the estimated signals for SNR = 20 dB.
The average squared angular error (SAE)  was used as the performance measure. Figure \ref{fig_ex_sine_exp} shows the simulation results, illustrating that:
\setlength{\leftmarginiii}{-2em}
\begin{itemize}
 \setlength{\itemindent}{-1em}
\item {\bf{PCA}} failed since the mixing vectors were not orthogonal and the source signals were correlated, both violating the assumptions for PCA.

\item {\bf{ICA}} (using the JADE algorithm \cite{cardoso1993blind}) failed because the signals were not statistically independent, as assumed in ICA.

\item {\bf{Low rank tensor approximation:}} a rank-2 CPD was used to estimate ${\bA}$ as the third factor matrix, which was then inverted to yield the sources. The accuracy of CPD was compromised as the components of tensor $\tX$ cannot be represented by rank-1 terms.

\item {\bf Low multilinear rank approximation:} Tucker decomposition (TKD) for the multilinear rank $(4,4,2)$ was able to retrieve the column space of the mixing matrix but could not find the individual mixing vectors due to the non-uniqueness of TKD.

\item {\bf BTD in multilinear rank-$(2,2,1)$ terms } matched the data structure \cite{LathauwerTBSS}, and it is remarkable that the sources were recovered using as few as 6 samples in the noise-free case.
\end{itemize}

%\begin{figure}[ht!]
%\centering
%%
%\begin{subfigure}[b]{.49\linewidth}
%      \centering
%{
%\includegraphics[width=1\linewidth, trim = 0cm 0.0cm 0cm 0cm , clip = true]{Fig6a}}
%%              \scriptsize   SAE$_{TD}$ = 4.0 dB, SAE$_{BCR Lo1}$ = 31.7 dB.
%%
%\end{subfigure}
%\hfill
%\begin{subfigure}[b]{.49\linewidth}
%      \centering
%{
%\includegraphics[width=1\linewidth, trim = 0cm 0.0cm 0cm 0cm , clip = true]
%{Fig6b}}
%\end{subfigure}
%\\[1ex]\hfil
%\begin{subfigure}[b]{.49\linewidth}
%     \centering
%{
%\includegraphics[width=1\linewidth, trim = 0cm 0.0cm 0cm 0cm , clip = true]
%{Fig6c}}
%\end{subfigure}
%%
%\hfill
%\begin{subfigure}[b]{.49\linewidth}
%      \centering
%{
%\includegraphics[width=1\linewidth, height =.67\linewidth,trim = 0cm 0.0cm 0cm 0cm , clip = true]{Fig6d}\label{fig_sine_exp_60samples_SAEmc}}
%%
%\end{subfigure}
%\captionsetup{singlelinecheck=off}
%\caption[foo bar]{Blind separation of the mixture of a pure sine wave and an exponentially modulated sine wave using PCA, ICA, CPD, Tucker decomposition (TKD)
%and BTD. The sources $s_1$ and $s_2$ are correlated and of  short duration; the symbols $\hat s_1$ and $\hat s_2$ denote the estimated sources.}
%\label{fig_ex_sine_exp}
%\end{figure}

\section{Higher-Order Compressed Sensing}

The aim of \textit{Compressed Sensing} (CS) is to provide  faithful reconstruction of a signal of interest when the set of available measurements is (much) smaller than the size of the original signal \cite{candes-2006,candes2006near,donoho2006compressed,Eldar:2012wf}. Formally, we have available $M$ (compressive) data samples $\mathbf{y} \in \Real^{M}$, which are assumed to be linear transformations of the original signal $\mathbf{x} \in \Real^{I}$ ($M<I$). In other words, $\mathbf{y}=\mathbf{\Phi} \mathbf{x}$, where the {\it sensing matrix} $\mathbf{\Phi}\in{\Real^{M\times I}}$ is usually random. Since the projections are of a lower dimension than the original data, the reconstruction is an ill-posed inverse problem, whose solution requires knowledge of the physics of the problem converted into constraints. For example, a 2D image $\bX \in \Real^{I_1\times I_2}$ can be vectorized as a long vector $\mathbf{x} = vec(\bX) \in \Real^{I}$ ($I=I_1I_2$) that admits  sparse representation in a known \textit{dictionary} $\mathbf{B}\in \Real^{I\times I}$,  so that $\mathbf{x} = \mathbf{B} \mathbf{g}$, where the matrix $\mathbf{B}$ may be a wavelet or discrete cosine transform (DCT) dictionary.
Then, faithful recovery of the original signal $\mathbf{x}$ requires finding the sparsest vector $\mathbf{g}$ such that:
\begin{equation}\label{CSproblem}
\mathbf{y}=\mathbf{W} \mathbf{g}, \mbox{  with } \|\mathbf{g}\|_0 \le K,  \qquad \mathbf{W}=\mathbf{\Phi}\mathbf{B},
\end{equation}
where $\|\cdot\|_0$ is the $\ell_0$-norm (number of non-zero entries) and $K\ll I$.

Since the $\ell_0$-norm minimization is not practical, alternative solutions involve iterative refinements of the estimates of vector $\mathbf{g}$ using greedy algorithms such as the Orthogonal Matching Pursuit (OMP) algorithm, or the \remove{use of}  $\ell_1$-norm minimization algorithms ($\|\mathbf{g}\|_1=\sum_{i=1}^I |g_i|$)  \cite{Eldar:2012wf}. Low coherence of the composite dictionary matrix $\mathbf{W}$ is a prerequisite for a satisfactory recovery of $\mathbf{g}$ (and hence $\mathbf{x}$) --- we need to choose $\mathbf{\Phi}$ and $\mathbf{B}$ so that the correlation between the columns of $\mathbf{W}$ is minimum \cite{Eldar:2012wf}.

When extending the CS framework to tensor data, we face two obstacles:
\begin{itemize}
\item {\it Loss of information}, such as spatial and contextual relationships in data, when a tensor $\tX \in \Real^{I_1\times I_2 \times \cdots \times I_N}$ is vectorized.
     %as $\mathbf{x} = vec(\tX) \in \Real^{I}$ where $I=I_1I_2\cdots I_N$; see also Section {\bf Tensorization}.
\item {\it Data handling}, since the size of  vectorized data and the associated dictionary $\mathbf{B}\in \Real^{I\times I}$ easily becomes  prohibitively large (see Section {\bf Curse of Dimensionality}), especially for tensors of high order.
\end{itemize}

\remove{Fortunately, tensor data are often highly structured, and thus particularly suited for compressive sampling. The benefits include reduced data acquisition efforts, compact storage and data completion (inpainting of entries that are missing, e.g., due to a broken sensor or unreliable measurement).}
\add{Fortunately, tensor data are typically highly structured -- a perfect match for compressive sampling -- so that the CS framework relaxes data acquisition requirements, enables compact storage, and facilitates data completion (inpainting of missing samples due to a broken sensor or unreliable measurement).}

\begin{figure}[t!]
\center
\includegraphics[width=12.4cm]{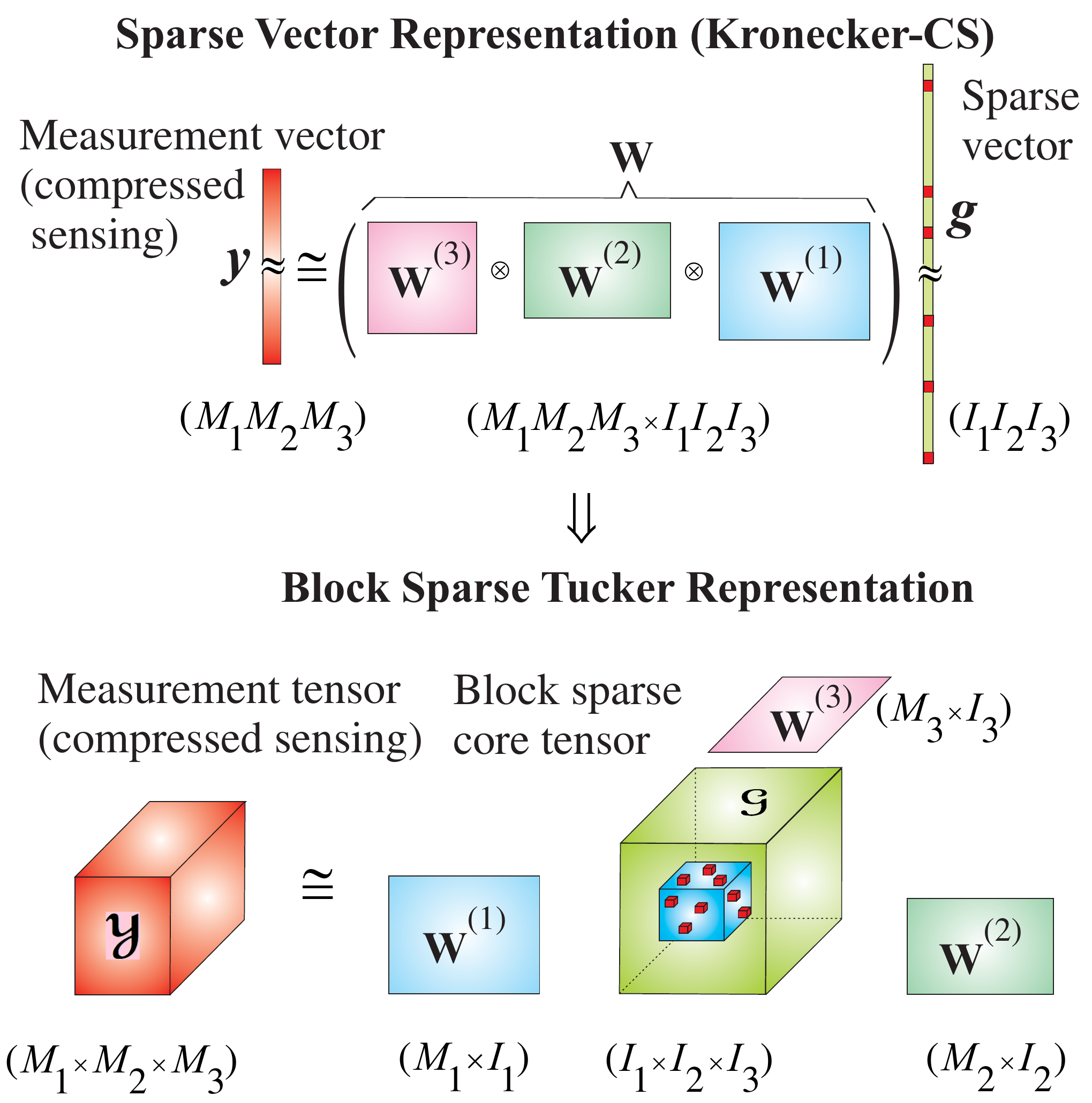}
\vspace{.2ex}
\caption{Compressed sensing with a Kronecker-structured dictionary. {\em Top:} Vector representation. {\em Bottom:} Tensor representation; Orthogonal Matching Pursuit (OMP) can perform  faster if the sparse entries belong to a small subtensor, up to permutation of the columns of $\bW^{(1)}$, $\bW^{(2)}$, $\bW^{(3)}$.}
\label{fig:Tuck-Kron-equiv}
\end{figure}

{\bf Kronecker-CS for fixed dictionaries.}
In many applications, the dictionary and the sensing matrix admit a Kronecker structure (Kronecker-CS model), as illustrated in \replace{Fig.}{Figure} \ref{fig:Tuck-Kron-equiv} (top) \cite{Duarte2012}. In this way, the global {\it composite dictionary matrix} becomes $\bW = \bW^{(N)} \otimes \bW^{(N-1)}\otimes \cdots \otimes \bW^{(1)}$, where each term $\mathbf{W}^{(n)}=\mathbf{\Phi}^{(n)}\mathbf{B}^{(n)}$ has a reduced dimensionality since $\mathbf{B}^{(n)}\in \Real^{I_n\times I_n}$ and $\mathbf{\Phi}^{(n)}\in \Real^{M_n\times I_n}$. Denote  $M=M_1M_2\cdots M_N$ and $I=I_1 I_2 \cdots I_N$, and \replace{knowing that}{since} $M_n\le I_n$, $n=1,2,\ldots,N$, this reduces storage requirements by a factor $\frac{\sum_n I_nM_n}{MI}$. The computation of $\mathbf{W} \mathbf{g}$ is affordable since $\mathbf{g}$ is sparse, however, computing $\mathbf{W}^T \mathbf{y}$ is expensive but can be efficiently implemented through a sequence of products involving much smaller matrices $\bW^{(n)}$ \cite{Caiafa2012-NC}. We refer to \cite{Duarte2012} for links between the coherence of factors $\bW^{(n)}$ and the coherence of the global composite dictionary matrix $\mathbf{W}$.

\replace{Fig.}{Figure} \ref{fig:Tuck-Kron-equiv} and Table \ref{CPD-Tucker} illustrate that the Kronecker-CS model is effectively a vectorized Tucker decomposition with a sparse core. The tensor equivalent of the CS paradigm in (\ref{CSproblem}) is therefore to find the sparsest core tensor $\tG$ such that:
\begin{equation}\label{TuckerCSproblem}
\tY \cong \tG \times_1 \bW^{(1)} \times_2 \bW^{(2)}  \cdots \times_N \bW^{(N)},
%\mbox{    with   } \|\tG\|_0 \le K,
\end{equation}
with $ \|\tG\|_0 \le K$, for a given set of mode-wise dictionaries $\bB^{(n)}$ and sensing matrices $\mathbf{\Phi}^{(n)}$ ($n=1,2,\dots,N$). Working with several small dictionary matrices, appearing in a Tucker representation, instead of a large global dictionary matrix, is an example of the use of tensor structure for efficient representation, see \add{also} Section {\bf Curse of Dimensionality}.

%
%\begin{figure}[t!]
%\center
%\def\svgwidth{1\linewidth}
%\input{Fig7.eps_tex}
%%\includegraphics[width=7.6cm]{fig7_r1}
%\vspace{.2ex}
%\caption{Compressed sensing with a Kronecker-structured dictionary. {\em Top:} Vector representation. {\em Bottom:} Tensor representation; Orthogonal Matching Pursuit (OMP) can perform  faster if the sparse entries belong to a small subtensor, up to permutation of the columns of $\bW^{(1)}$, $\bW^{(2)}$, $\bW^{(3)}$.}
%\label{fig:Tuck-Kron-equiv}
%\end{figure}

A higher-order extension of the OMP algorithm, referred to as the Kronecker-OMP algorithm \cite{Caiafa2012-NC}, requires $K$ iterations to find the $K$ non-zero entries of the core tensor $\tG$. Additional computational advantages can be gained if it can be assumed that the $K$ non-zero entries belong to a small subtensor of $\tG$,
as shown in \replace{Fig.}{Figure} \ref{fig:Tuck-Kron-equiv} (bottom); such a structure is inherent to e.g., hyperspectral imaging \cite{Caiafa2012-NC, Caiafa:2013} and 3D astrophysical signals. More precisely, if the $K=L^N$ non-zero entries are located within a subtensor of size ($L\times L \times \cdots \times L$), where $L \ll I_n$, then the so-called $N$-way Block OMP algorithm (N-BOMP) requires at most $NL$  iterations, which is linear in $N$ \cite{Caiafa2012-NC}. The Kronecker-CS model has been applied in Magnetic Resonance Imaging (MRI), hyper-spectral imaging, and in the inpainting of multiway data \cite{Caiafa:2013,Duarte2012}.

{\bf Approaches without fixed dictionaries.} In Kronecker-CS the mode-wise dictionaries $\mathbf{B}^{(n)}\in \Real^{I_n\times I_n}$ can be chosen so as best to represent physical properties or prior knowledge about the data. They can also be learned from a large ensemble of data tensors, for instance in an ALS type fashion \cite{Caiafa:2013}. Instead of the total number of sparse entries in the core tensor, the size of the core (i.e., the multilinear rank) may be used as a measure for sparsity so as to obtain a low-complexity representation from compressively sampled data \cite{Gandy2011,signoretto2013learning}. Alternatively, a PD representation can be used instead of a Tucker representation. Indeed, early work in chemometrics involved excitation-emission data for which part of the entries was unreliable because of scattering; the CPD of the data tensor is then computed by treating such entries as missing \cite{Smilde}. While CS variants of several CPD algorithms exist \cite{Acar-Morup11,Sorber-tensorlab}, the ``oracle'' properties of tensor-based models are still not as well understood as for their standard models; a notable exception is CPD with sparse factors \cite{Kyrillidis2012}.
\begin{figure}[ht!]
\center
\includegraphics[width=11.99cm]{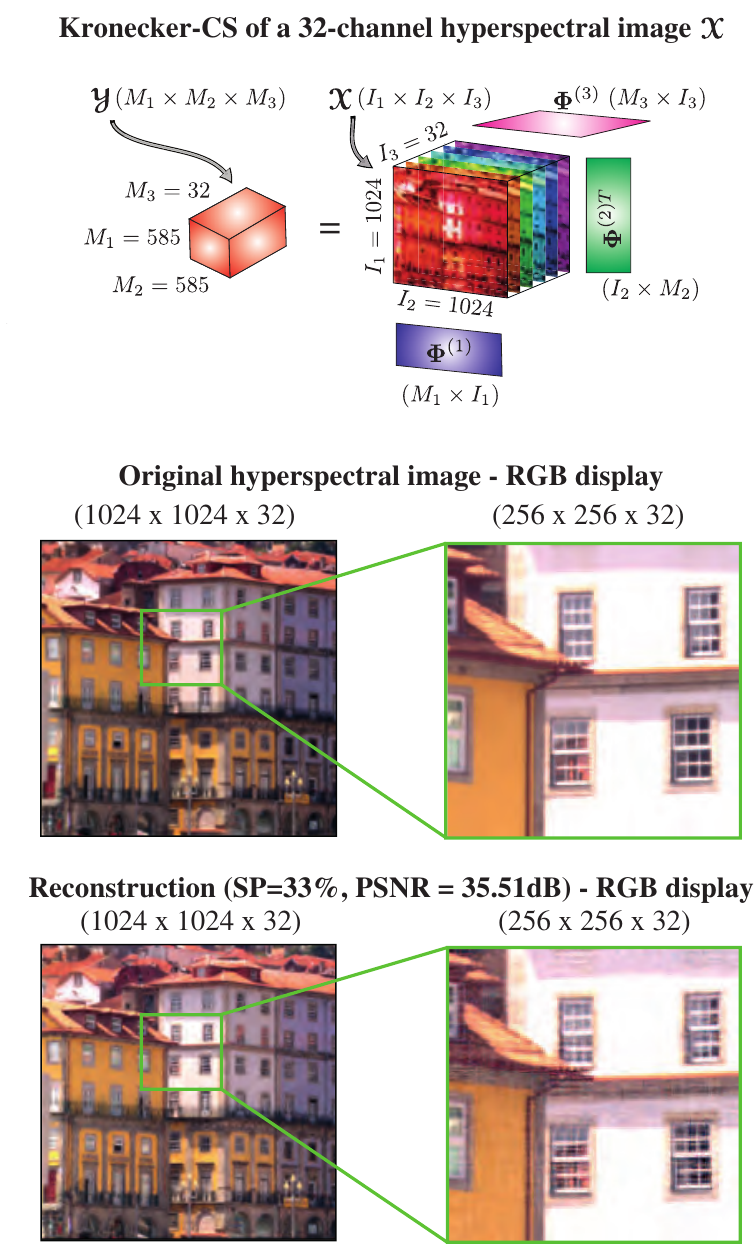}
\caption{Multidimensional compressed sensing of a 3D hyperspectral image using Tucker representation with a small sparse core in wavelet bases.}
\label{fig:CSexample}
\end{figure}

{\bf Example 2.} \replace{Fig.}{Figure} \ref{fig:CSexample} shows an original 3D ($1024 \times 1024 \times 32$) hyperspectral image $\tX$ which contains scene reflectance measured at 32 different frequency channels, acquired by a low-noise Peltier-cooled digital camera in the wavelength range of 400--720 nm \cite{Foster-DB}. Within the Kronecker-CS setting, the tensor of compressive measurements $\tY$ was obtained by multiplying the frontal slices by random Gaussian sensing matrices $\mathbf{\Phi}^{(1)}\in{\Real^{M_1\times 1024}}$ and $\mathbf{\Phi}^{(2)} \in{\Real^{M_2\times 1024}}$ ($M_1, M_2<1024$) in the first and second mode, respectively, while $\mathbf{\Phi}^{(3)}\in{\Real^{32\times 32}}$ was the identity matrix (see \replace{Fig.}{Figure}~\ref{fig:CSexample} (top)). We used  Daubechies wavelet factor matrices $\mathbf{B}^{(1)}=\mathbf{B}^{(2)}\in{\Real^{1024\times 1024}}$ and $\mathbf{B}^{(3)}\in{\Real^{32\times 32}}$, and \replace{used}{employed} N-BOMP to recover the small sparse core tensor and, subsequently, reconstruct the original 3D image as shown in \replace{Fig.}{Figure}~\ref{fig:CSexample} (bottom). For the Sampling Ratio SP=$33\%$ ($M_1=M_2=585$) this gave the Peak Signal to Noise Ratio (PSNR) of $35.51$dB, while taking $71$ minutes to compute the required $N_{iter}=841$ sparse entries. For the same quality of reconstruction (PSNR=$35.51$dB), the more conventional Kronecker-OMP algorithm found $0.1$\% of the wavelet coefficients as significant, thus requiring $N_{iter}=K=0.001\times (1024\times 1024 \times 32) = 33,555$ iterations and days of computation time.

\section{Large-Scale Data and Curse of Dimensionality}

The sheer size of tensor data easily exceeds the memory or saturates the processing capability of standard  computers, it is therefore natural to ask ourselves how tensor decompositions can be computed if the tensor dimensions in all or some modes are large or, worse still, if the tensor order is high. The term \textit{curse of dimensionality}, in a general sense, was introduced by %Richard E.
Bellman to refer to various computational bottlenecks when dealing with high-dimensional settings\remove{, more specifically,}. In the context of tensors, the \textit{curse of dimensionality} refers to the fact that the number of elements of an $N$th-order $(I \times I \times \cdots \times I)$ tensor, $I^N$, scales {\emph{exponentially}} with the tensor order $N$. For example, the number of values of a discretized function in Figure \ref{Fig:Tensorization} (bottom),  quickly becomes unmanageable in terms of both computations and storing as $N$ increases. In addition to their standard use (signal separation, enhancement, etc.), tensor decompositions may be elegantly employed in this context as efficient representation tools. The first question is then which type of tensor decomposition is appropriate.

{\bf Efficient data handling.} If all computations are performed on a {\em CP representation} and not on the raw data tensor itself, then instead of the original $I^N$  raw data entries, the number of parameters in a CP representation reduces to $N I R$, which scales {\emph{linearly}} in $N$ (see Table \ref{table_complexity}). This effectively bypasses the \textit{curse of dimensionality},
while  giving us the freedom to choose the rank $R$ as a function of the desired accuracy \cite{beylkin}; on the other hand the CP approximation may involve numerical problems (see Section {\bf Canonical Polyadic Decomposition}).

Compression is also inherent to {\em Tucker decomposition}, as it reduces the size of a given data tensor from the original $I^N$ to $(NIR+R^N)$, thus exhibiting an approximate compression ratio of $(\frac{I}{R})^N$. We can then benefit from the well understood and reliable approximation by means of matrix SVD, however, this is only meaningful for low $N$.

\minrowclearance 2ex
\begin{table}[ht]
\centering
\caption{Storage complexities of  tensor models for an $N$th-order  tensor
$\tX \in \Real^{I \times I \times \cdots \times I}$, whose original storage complexity is ${\cal{O}}(I^N)$.}
%\centering
%
{\small \shadingbox{
\begin{tabular*}{0.4\linewidth}[t]{@{\extracolsep{\fill}}ll} \hline
1.  CPD   &  ${\cal{O}}(NIR)$ \\
2. Tucker   & ${\cal{O}}(NIR +R^N)$ \\
3. TT   &   ${\cal{O}}(NIR^2)$ \\
4. QTT   & ${\cal{O}}(N R^2 \log_2(I))$ \\
\hline
    \end{tabular*}
  }}
%    \end{center}
\label{table_complexity}
\end{table}
\minrowclearance 0ex

{\bf Tensor networks.} A numerically reliable way to tackle curse of dimensionality is through a concept from scientific computing and quantum  information theory, termed {\em tensor networks}, which  represents a tensor of a possibly very high order as a set of sparsely interconnected matrices and core tensors of low order (typically, order 3). %\cite{Cichocki-SISA}.
% \cite{Cichocki-SISA,Orus2013}(see also Table \ref{table_notation2}).
These low-dimensional  cores are interconnected via tensor contractions to provide a highly compressed representation of a data tensor. In addition, existing algorithms for the approximation of a given tensor by a tensor network have good numerical properties, making it possible to control the error and achieve any desired accuracy of approximation.
For example, tensor networks allow for the representation of a wide class of discretized multivariate functions even in cases where the number of function values is larger than  the number of atoms in the universe \cite{Khoromskij-TT,Hackbush2012,Grasedyck-rev}.
\begin{figure}[t]
\centering
\includegraphics[width=0.6\linewidth]{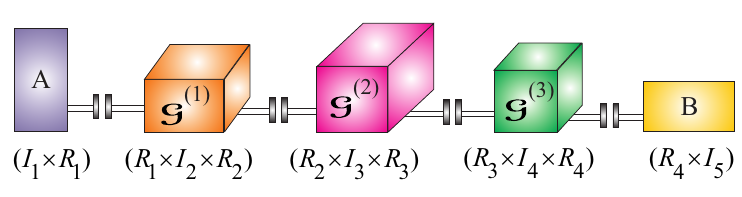}\\
\caption{Tensor Train (TT) decomposition of a fifth-order  tensor $\tX \in \Real^{I_1 \times I_2 \times \cdots \times I_5}$, consisting of two matrix carriages and three third-order tensor carriages. The five carriages are connected through tensor contractions, which can be expressed in a scalar form as $x_{i_1,i_2,i_3,i_4,i_5} =\sum_{r_1=1}^{R_1} \sum_{r_2=1}^{R_2} \cdots \sum_{r_4=1}^{R_4} a_{i_1,r_1} \; g^{(1)}_{r_1,i_2,r_2} \; g^{(2)}_{r_2,i_3,r_3} \; g^{(3)}_{r_3,i_4,r_4} \; b_{r_4,i_5} $.}
%The $(N-1)$-tuple $(R_1,R_2,\ldots,R_{N-1})$ is  called the TT-rank.}}
\label{Fig:TT5}
\end{figure}

Examples of tensor networks are the hierarchical Tucker (HT) decompositions and Tensor Trains (TT) (see Figure \ref{Fig:TT5}) \cite{hTucker,OseledetsTT11}. The TTs are also known
as Matrix Product States (MPS) and have been used by physicists more than  two decades  (see \cite{Cichocki-SISA,Orus2013} and references therein).
The PARATREE algorithm was developed in signal processing and follows a similar idea, it uses a polyadic representation of a data tensor (in a possibly nonminimal  number of terms), whose  computation then requires only the matrix SVD \cite{VisaSP-09}.

\begin{figure*}[t]
\centering
\includegraphics[width=.95\textwidth]{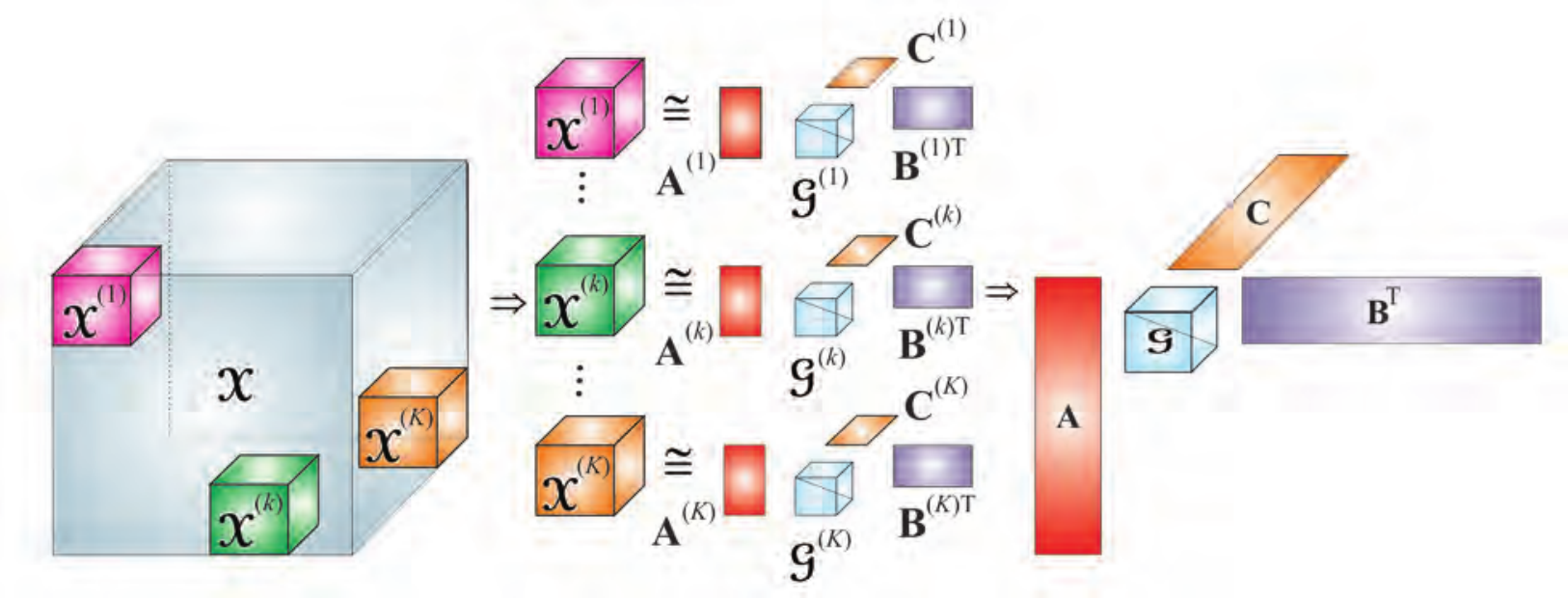}
\caption{Efficient computation of the CP and Tucker decompositions, whereby tensor decompositions are computed in parallel for sampled blocks, \replace{which}{these} are then merged to obtain the global components $\bA$, $\bB$, $\bC$ and a core tensor $\tG$.}
\label{Fig:blocks}
\end{figure*}

%\begin{figure*}[ht!]
%\centering
%\def\svgwidth{1\linewidth}
%\input{Fig10.eps_tex}
%%\includegraphics[width=.8\textwidth]{fig10_r1c}
%\caption{Efficient computation of the CP and Tucker decompositions, whereby tensor decompositions are computed in parallel for sampled blocks, \replace{which}{these} are then merged to obtain the global components $\bA$, $\bB$, $\bC$ and a core tensor $\tG$.}
%\label{Fig:blocks}
%\end{figure*}

For very large-scale data that exhibit a well-defined structure, an even more radical approach can be employed to achieve a parsimonious representation --- through the concept of {\em quantized} or {\em quantic tensor networks} (QTN)  \cite{Khoromskij-TT,Grasedyck-rev}. For example, a huge vector $\bx \in \Real^I$ with $I = 2^L$ elements can be quantized and tensorized through reshaping into a $(2 \times 2  \times \cdots \times 2)$  tensor $\tX$ of order $L$, as illustrated in Figure \ref{Fig:Tensorization} (top). If $\bx$ is an exponential signal, $x(k) = a z^{k}$, then $\tX$ is a symmetric rank-1 tensor that can be represented by two parameters: the scaling factor $a$ and the generator $z$ (cf. (\ref{Eq:Hankelization_Exponential}) in Section {\bf Tensorization}). Non-symmetric terms provide further opportunities, beyond the sum-of-exponential representation by symmetric low-rank tensors. Huge matrices and tensors may be dealt with in the same manner.
For instance,  an $N$th-order tensor $\tX \in \Real^{I_1 \times \cdots \times I_N}$, with  $I_n=q^{L_n}$,  can be quantized  in all modes simultaneously  to yield a
$(q \times q \times \cdots \times q)$  quantized tensor of higher order. In QTN, $q$ is small, typically $q=2,3,4$,
for example, the binary encoding ($q=2$) reshapes an $N$th-order tensor with $(2^{L_1} \times 2^{L_2} \times \cdots \times 2^{L_N})$ elements into a tensor of order $(L_1+L_2 +\cdots +L_N)$ with the same number of elements.
The tensor train decomposition applied to quantized tensors is
referred to as the \add{quantized TT (QTT)}; variants for other tensor representations have also been derived \cite{Khoromskij-TT,Grasedyck-rev}. In scientific computing, such formats provide the so-called \emph{super-compression} --- a logarithmic reduction of storage requirements: ${\cal{O}}(I^N) \rightarrow {\cal{O}}(N \log_q(I))$.

{\bf Computation of the decomposition/representation.} Now that we have addressed the possibilities for efficient tensor representation, the question that needs to be answered is how these representations can be computed from the data in an efficient manner. The first approach is to process the data in smaller blocks rather than in a batch manner \cite{Phan-CP}. In such a ``divide-and-conquer'' approach, different blocks may be processed in parallel and their decompositions carefully recombined (see \replace{Fig.}{Figure} \ref{Fig:blocks}) \cite{Suter13,Phan-CP}. In fact, we may even compute the decomposition through recursive updating, as new data arrive \cite{sid_adPARAFAC}. Such recursive techniques may be used for efficient computation and for tracking decompositions in the case of nonstationary data.

The second approach would be to employ compressed sensing ideas (see Section {\bf Higher-Order Compressed Sensing}) to fit an algebraic model with a limited number of parameters to possibly large data. In addition to completion, the goal here is a significant reduction of the cost of data acquisition, manipulation and storage --- breaking the Curse of Dimensionality being an extreme case.

\begin{figure}[t!]
\centering
\includegraphics[width=9.6cm]{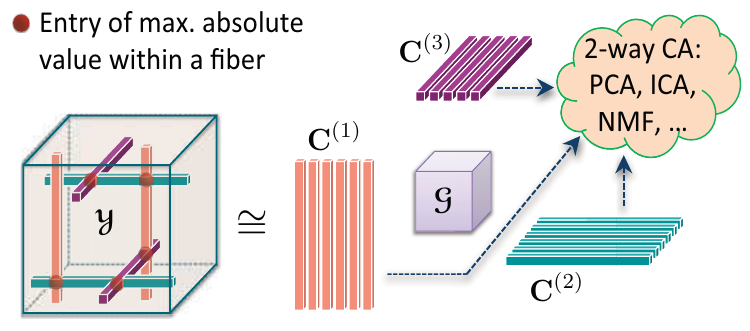}
\caption{Tucker representation through fiber sampling and cross-approximation: the columns of  factor matrices are sampled from the fibers of the original data tensor $\tX$. Within MWCA the selected fibers may be further processed using BSS algorithms.
}
\label{fig:MWCA}
\end{figure}

While algorithms for this purpose are available both for low rank and low multilinear rank representation \cite{Acar-Morup11,Gandy2011}, an even more drastic approach would be to directly adopt sampled fibers as the bases in a tensor representation. In the Tucker decomposition setting we would choose the columns of the factor matrices $\bB^{(n)}$ as mode-$n$ fibers of the tensor, which requires addressing the following two problems: (i) how to find fibers that allow us to best represent the tensor, and (ii) how to compute the corresponding core tensor at a low cost (i.e., with minimal access to the data).
The matrix counterpart of this problem (i.e., representation of a large matrix on the basis of a few columns and rows) is referred to as the {\emph{pseudoskeleton approximation}} \cite{Goreinov:1997}, where the optimal representation corresponds to the columns and rows that intersect in the submatrix of maximal volume (maximal absolute value of the determinant). Finding the optimal submatrix is computationally hard, but quasi-optimal submatrices may be found by heuristic so-called ``cross-approximation'' methods that only require a limited, partial exploration of the data matrix. Tucker variants of this approach have been derived in \cite{Caiafa-Cichocki-CUR,Goreinov:2008hd,Oseledets2008} and are illustrated in \replace{Fig.}{Figure}~\ref{fig:MWCA}, while cross-approximation for the TT format has been derived in \cite{oseledets2010tt}.  Following a somewhat different idea, a tensor generalization of the CUR decomposition of matrices samples fibers on the basis of statistics derived from the data \cite{mahoney2008tensor}.

\section{Multiway Regression --- Higher Order PLS (HOPLS)}
\label{sec:HOPLS}
{\bf Multivariate regression.} Regression refers to the modelling of one or more {\em dependent variables (responses)}, $Y$, by a set of {\em independent data (predictors)}, $X$. In the simplest case of conditional MSE estimation,  ${\hat y}=E(y|x)$, the response $y$ is a linear combination of the elements of the vector of predictors ${\bf x}$; for multivariate data the Multivariate Linear Regression (MLR) uses a matrix  model, ${\bf Y}={\bf X} {\bf P} + {\bf E}$, where ${\bf P}$ is the matrix of coefficients (loadings) and ${\bf E}$ the residual matrix. The  MLR solution gives ${\bf P} = \big( {\bf X}^{T} {\bf X} \big)^{-1}{\bf X}^{T} {\bf Y}$, and involves  inversion of the moment matrix ${\bf X}^{T} {\bf X}$.
A common technique to stabilize the inverse of the moment matrix ${\bf X}^{T} {\bf X}$ is {\em principal component regression (PCR)}, which employs low rank approximation of ${\bf X}$.

\begin{figure}[t!]
\centering
\includegraphics[width=9.6cm]{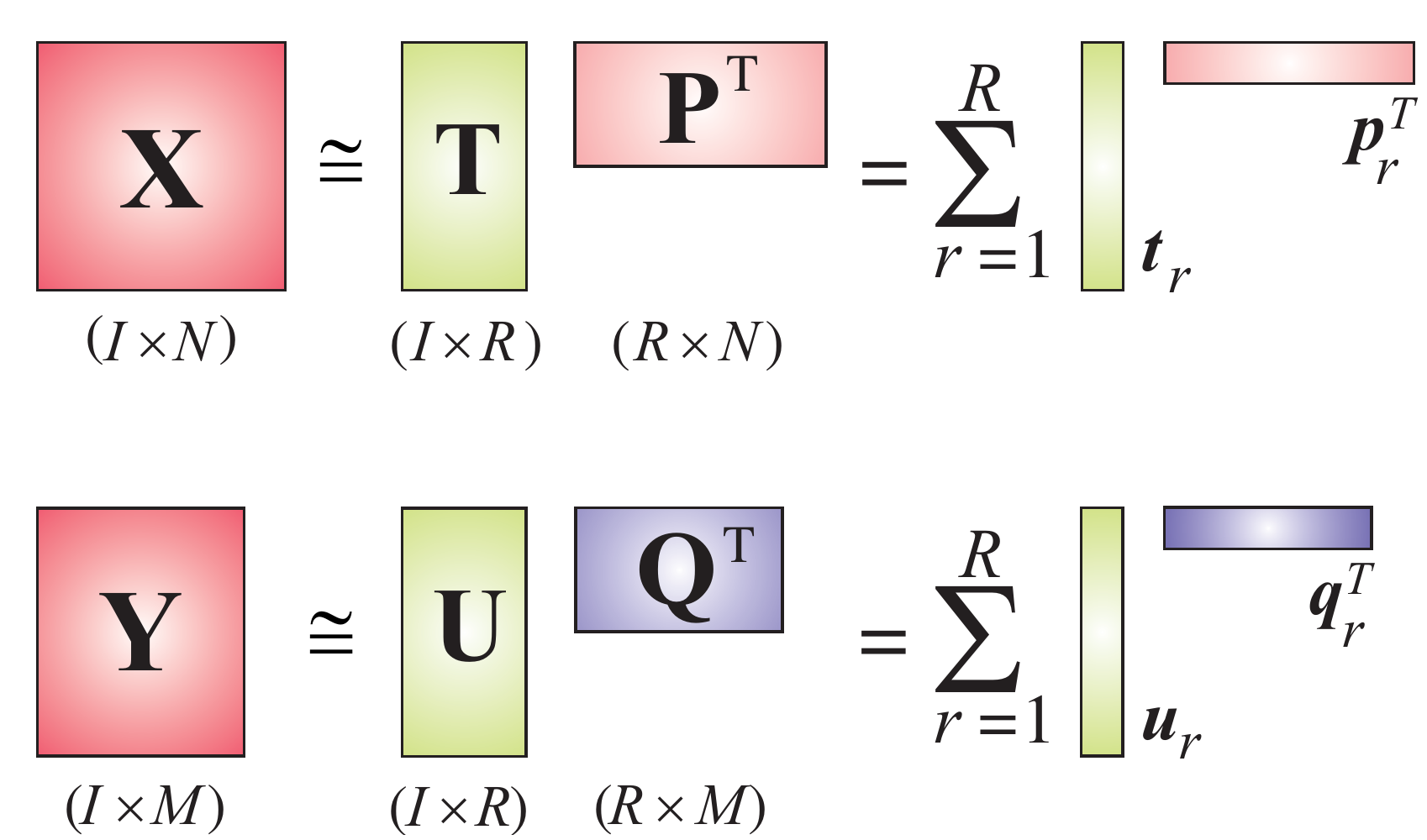}
\caption{The basic PLS model performs joint sequential low-rank approximation of the matrix of predictors $\bX$ and the matrix of responses $\bY$, so as to share (up to the scaling ambiguity) the latent components --- columns of the {\em score matrices} $\bT$ and $\bU$.   The matrices $\bP$ and $\bQ$ are the {\em loading matrices} for  predictors and responses, and $\bE$ and $\bF$ are the corresponding {\em residual matrices}. }
\label{fig:PLS}
\end{figure}

{\bf Modelling structure in data --- the PLS.} Notice that in stabilizing multivariate regression PCR uses only  information in the $X$-variables, with no feedback from the $Y$-variables.
%  \remove{and thus no means to reflect the structure in data or guarantee parsimonious representation and interpretability of the components.}
The idea behind the Partial Least Squares (PLS) method  is to account for structure in data by assuming that the underlying system is governed by a small number, $R$, of specifically constructed latent variables, called {\em scores}, that are shared between the $X$- and $Y$-variables; in estimating the number $R$, PLS compromises between fitting ${\bf X}$ and predicting ${\bf Y}$. Figure {\ref{fig:PLS}} illustrates that the PLS procedure:  (i) uses eigenanalysis to perform {\em contraction} of the data matrix ${\bf X}$ to the principal eigenvector {\em score matrix} ${\bf T} =[{\bf t}_{1},\ldots, {\bf t}_{R}]$ of  rank $R$;
(ii) ensures that the ${\bf t}_{r}$ components are maximally correlated with the  ${\bf u}_{r}$ components in the approximation of the responses ${\bf Y}$, this is achieved when the ${\bf u}_{r}$'s are scaled versions of the ${\bf t}_{r}$'s. The $Y$-variables are then regressed on the matrix ${\bf U}=[{\bf u}_{1},\ldots,{\bf u}_{R}]$.
Therefore, PLS is a  multivariate model with inferential ability that
aims to find a representation of ${\bf X}$ (or a part of ${\bf X}$) that is relevant for predicting ${\bf Y}$, using the model
\begin{eqnarray}
\bX &=& \bT \; \bP^T + \bE = \sum_{r=1}^R \bt_r \; \bp_r^T +\bE, \label{eq:PLSX} \\
\bY &=& \bU \; \bQ^T + \bF = \sum_{r=1}^R  \bu_r \; \bq_r^T  + \bF. \label{eq:PLSY}
\end{eqnarray}
The score vectors ${\bf t}_{r}$  provide an LS fit of $\bX$-data, while at the same time  the maximum  correlation between ${\bf t}$- and ${\bf u}$-scores ensures a good predictive model for Y-variables.
The predicted responses $\bY_{new}$ are then obtained from  new data $\bX_{new}$
 and the loadings $\bP$ and  $\bQ$.

In practice, the score vectors ${\bf t}_{r}$ are extracted sequentially, by a series of orthogonal projections followed by the deflation of ${\bf X}$. Since the rank of ${\bf Y}$ is not necessarily decreased  with each new $\bt_{r}$, we may continue deflating until the rank of the ${\bf X}$-block is exhausted, so as to balance between prediction accuracy and model order.

%\begin{figure}[t!]
%\centering
%\def\svgwidth{1\linewidth}
%\input{fig12.eps_tex}
%%\includegraphics[width=7.0cm]{fig12}
%\caption{the basic pls model performs joint sequential low-rank approximation of the matrix of predictors $\bx$ and the matrix of responses $\by$, so as to share (up to the scaling ambiguity) the latent components --- columns of the {\em score matrices} $\bt$ and $\bu$.   the matrices $\bp$ and $\bq$ are the {\em loading matrices} for  predictors and responses, and $\be$ and $\bf$ are the corresponding {\em residual matrices}. }
%\label{fig:pls}
%\end{figure}

\begin{figure}[t!]
\centering
\includegraphics[width=12.8cm]{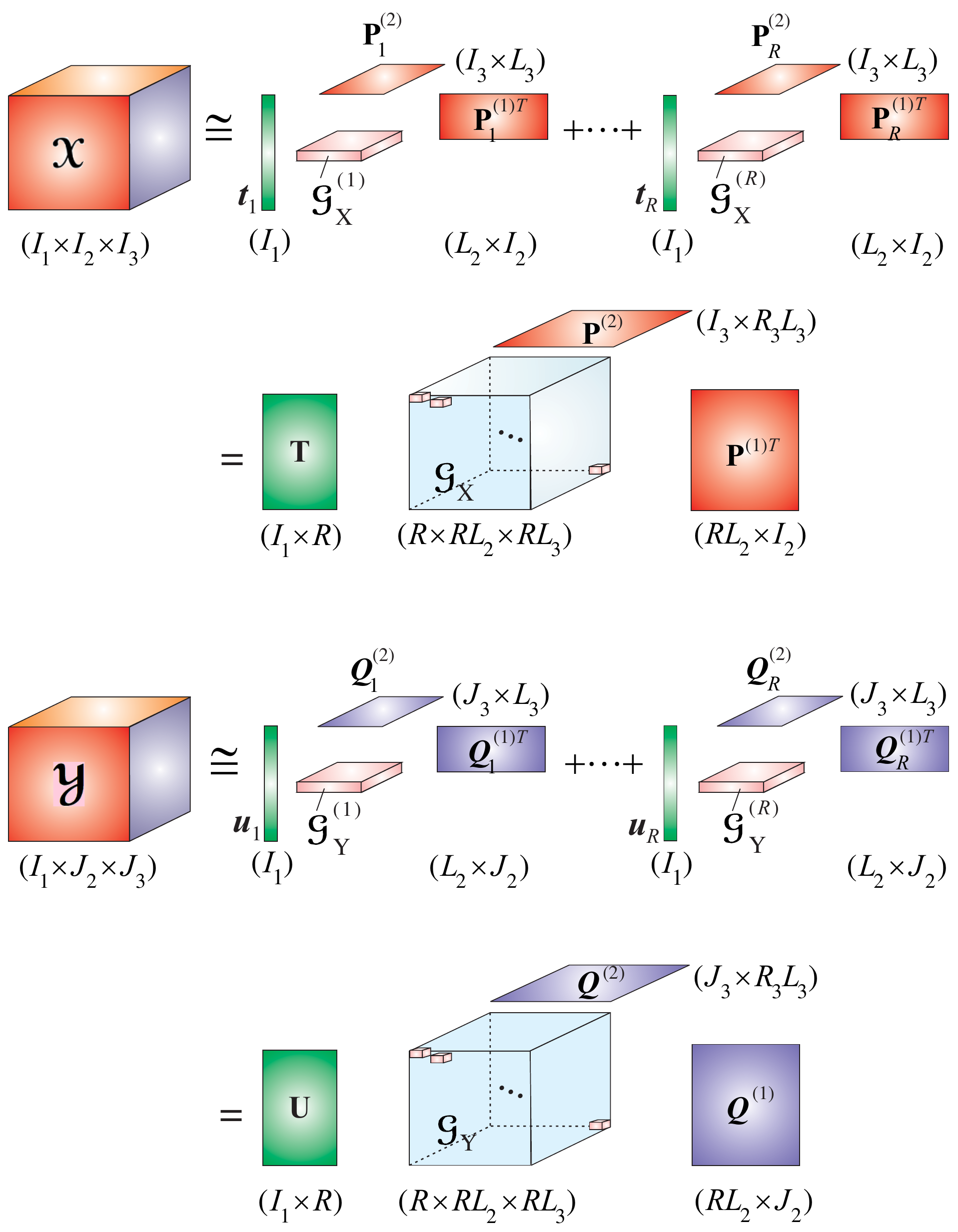}
\caption{The principle of Higher Order PLS (HOPLS) for third-order tensors. The core tensors $\tG_{\bX}$ and $\tG_{\bY}$ are block-diagonal.
The BTD-type structure allows for the modelling of general components that are highly correlated in the first mode.}
\label{fig:HOPLS1}
\end{figure}

The PLS concept can be generalized to tensors in the following ways:
\begin{enumerate}
\item {\em By unfolding  multiway data}. For example $\tX(I \times J \times K)$ and $\tY(I \times M \times N)$ can be flattened into long matrices ${\bf X}(I \times JK)$ and ${\bf Y}(I \times MN)$, so as to admit matrix-PLS  (see Figure {\ref{fig:PLS}}).  However, the flattening prior to standard {\em bilinear} PLS obscures structure in multiway data and compromises the interpretation of latent components.
\item {\em By low rank tensor approximation.} The so-called N-PLS attempts to find score vectors having maximal covariance with response variables, under the constraints that tensors $\tX$ and $\tY$ are decomposed as a sum of rank-one tensors \cite{Bro1996}.
\item {\em By a BTD-type approximation,} as in the Higher Order PLS  (HOPLS) model shown in Figure {\ref{fig:HOPLS1}} \cite{Qibin-HOPLS}. The use of block terms within HOPLS equips it with additional flexibility,
together with a more realistic analysis than unfolding-PLS and N-PLS.
\end{enumerate}
The principle of HOPLS can be formalized  as a set of sequential approximate decompositions of the independent tensor $\tX \in \Real^{I_1 \times I_2 \times \cdots \times I_N}$
and the dependent tensor $\tY \in \Real^{J_1 \times J_2 \times \cdots \times J_M}$ (with $I_1 = J_1$),
so as to ensure maximum similarity (correlation) between the scores ${\bf t}_{r}$ and ${\bf u}_{r}$ within the loadings matrices $\bT$ and $\bU$, based on
\begin{align}
\tX &\cong \sum_{r=1}^R \tG_{\bX}^{(r)}  \times_1 \bt_r \times_2 \bP^{(1)}_r \cdots  \times_N \bP^{(N-1)}_r  \quad \\ %\notag \\
\tY &\cong \sum_{r=1}^R \tG_{\bY}^{(r)}  \times_1 \bu_r  \times_2 \bQ^{(1)}_r  \cdots  \times_N \bQ^{(M-1)}_r.    %\\ %\notag
\end{align}
%%
%\begin{eqnarray}
%\tX &=& \sum_{r=1}^R \tG_{\bX}^{(r)}  \times_1 \bt_r \times_2 \bP^{(1)}_r \times_3 \cdots  \times_N \bP^{(N-1)}_r + \tE \\ %\notag \\
%\tY &=& \sum_{r=1}^R \tG_{\bY}^{(r)}  \times_1 \bu_r  \times_2 \bQ^{(1)}_r \times_3 \cdots  \times_N \bQ^{(M-1)}_r + \tF %\\ %\notag
%\end{eqnarray}
%%
%where the tensors $\tE$ and $\tF$ represent errors (residuals). Large residuals in the response, $\tF$, indicate a poor HOPLS model, while the predictor residuals, $\tE$, reflect the part of $\tX$ that is not used in the modelling of $\tY$.
A number of data-analytic problems can be reformulated as
either regression or  ``similarity analysis'' (ANOVA, ARMA, LDA, CCA), so that both the matrix and tensor PLS solutions can be generalized across  exploratory data analysis.

%\begin{figure}[t!]
%\centering
%\def\svgwidth{1\linewidth}
%\input{Fig13.eps_tex}
%%\includegraphics[width=\linewidth]{fig13_r1}
%\caption{The principle of Higher Order PLS (HOPLS) for third-order tensors. The core tensors $\tG_{\bX}$ and $\tG_{\bY}$ are block-diagonal.
%The BTD-type structure allows for the modelling of general components that are highly correlated in the first mode.}
%\label{fig:HOPLS1}
%\end{figure}

\begin{figure}[t]
\centering
%centerline{
%\def\svgwidth{1\linewidth}
%\input{fig14_r2b.pdf_tex}
%\input{fig14_r1.eps_tex}
\includegraphics[width=16.2cm]{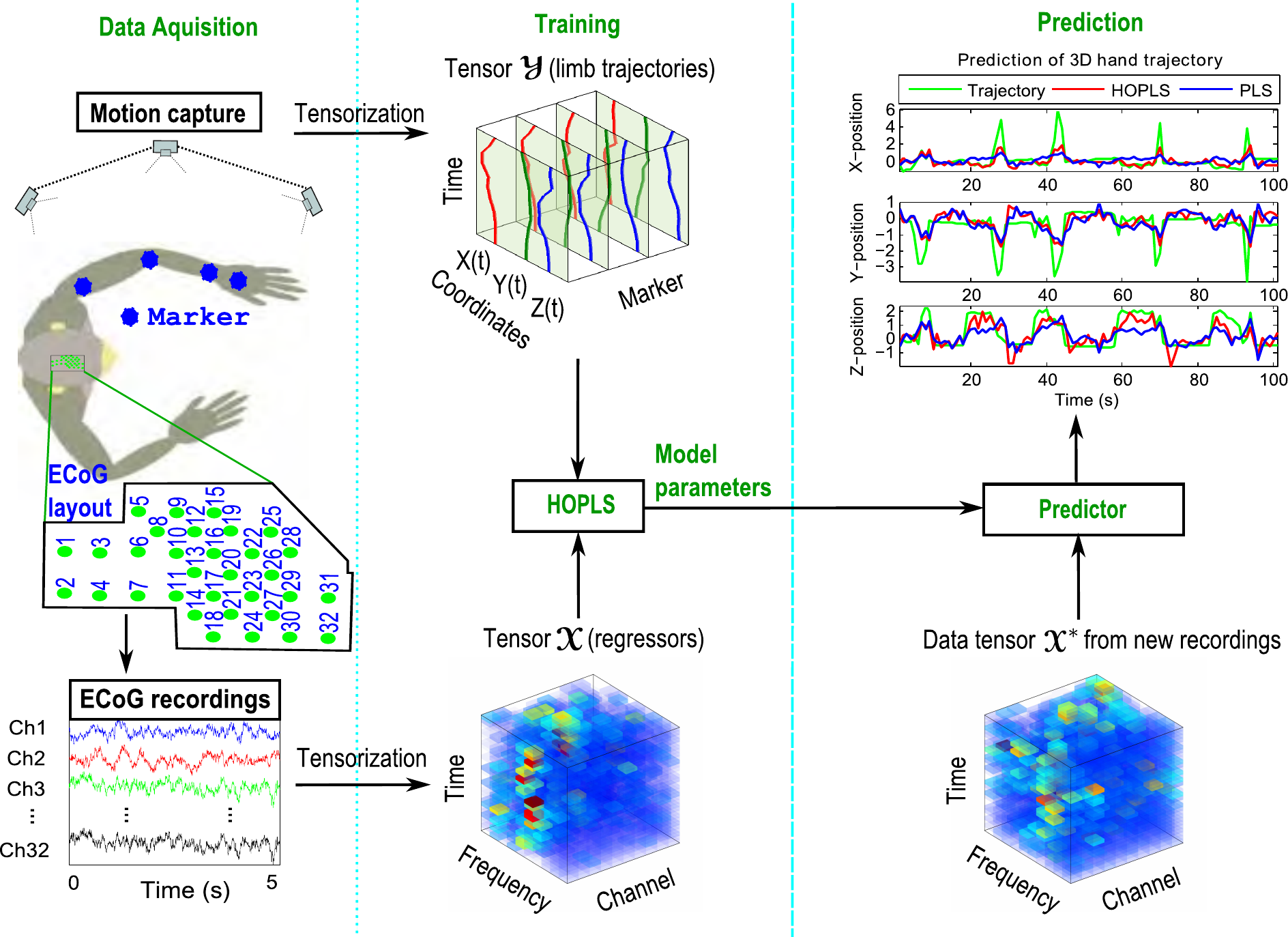}
  \caption{Prediction of arm movement from brain electrical responses. {\em Left:} Experiment setup. {\em Middle:} Construction of the data and response tensors and training. {\em Right:} The new data tensor (bottom) and the predicted 3D arm movement trajectories ($X, \;Y, \;Z$ coordinates) obtained by tensor-based HOPLS and standard matrix-based PLS (top).}
  \label{fig:ECoGdecoding}
\end{figure}

{\bf Example 4: Decoding of a 3D hand movement trajectory from the
electrocorticogram (ECoG).} The predictive power of tensor-based PLS is illustrated on a real-world example of the prediction of  arm movement trajectory from ECoG.  Fig\@. {\ref{fig:ECoGdecoding}}(left) illustrates the experimental setup, whereby 3D arm movement of a monkey was captured by an optical motion capture system with reflective markers affixed to the left shoulder, elbow,
wrist, and hand; for full detail see (\url{http://neurotycho.org}). The predictors (32 ECoG channels)   naturally build a fourth-order tensor $\tX$ (time$\times$channel\_no$\times$epoch\_length$\times$frequency) while the movement trajectories for the four markers (response) can be  represented as a third-order tensor $\tY$ (time$\times$3D\_marker\_position$\times$marker\_no). The goal of the training stage is to identify the HOPLS parameters: $\tG_{\bX}^{(r)}, \tG_{\bY}^{(r)}, {\bf P}^{(n)}_r,{\bf Q}^{(n)}_r$, see also Figure {\ref{fig:HOPLS1}}.  In the test stage, the movement trajectories, $\tY^*$, for the new ECoG data, $\tX^*$, are predicted through
multilinear projections: (i) the new scores, ${\bf t}_{r}^{*}$, are found from new data, $\tX^{*}$, and the existing model parameters: $\tG_{{\bf X}}^{(r)}, {\bf P}^{(1)}_r, {\bf P}^{(2)}_r, {\bf P}^{(3)}_r$,  (ii) the predicted trajectory is calculated as
$
\tY^* \approx \sum_{r=1}^R \tG_{\bY}^{(r)}  \times_1 \bt^*_r  \times_2 \bQ^{(1)}_r \times_3  \bQ^{(2)}_r \times_4 \bQ^{(3)}_r$.  In the simulations, standard PLS was applied in the same way to the unfolded tensors.

Figure~{\ref{fig:ECoGdecoding}}(right) shows that although the standard PLS was able to predict the movement corresponding to
each marker individually, such prediction is quite crude as the two-way PLS does not adequately account for mutual information among the four markers.
The enhanced predictive performance of the BTD-based HOPLS (red line in Fig\@.{\ref{fig:ECoGdecoding}}(right)) is therefore attributed to its ability to model interactions between complex latent components of both predictors and responses.

\section{Linked Multiway Component Analysis and Tensor Data Fusion}

Data fusion concerns joint analysis of an ensemble of data sets, such as multiple ``views'' of a particular phenomenon, where some parts of the ``scene'' may be visible in only one or a few data sets.
Examples include fusion of visual and thermal images in low visibility conditions, or the analysis of human electrophysiological signals in response to a certain stimulus but
from different subjects and trials; these are naturally analyzed together by means of matrix/tensor factorizations.
The ``coupled'' nature of the analysis of multiple datasets ensures that there may be common factors across the datasets, and that some components are not shared (e.g., processes that are independent of excitations or stimuli/tasks).

\begin{figure}[t!]
 \centering
\includegraphics[width=9.6 cm]{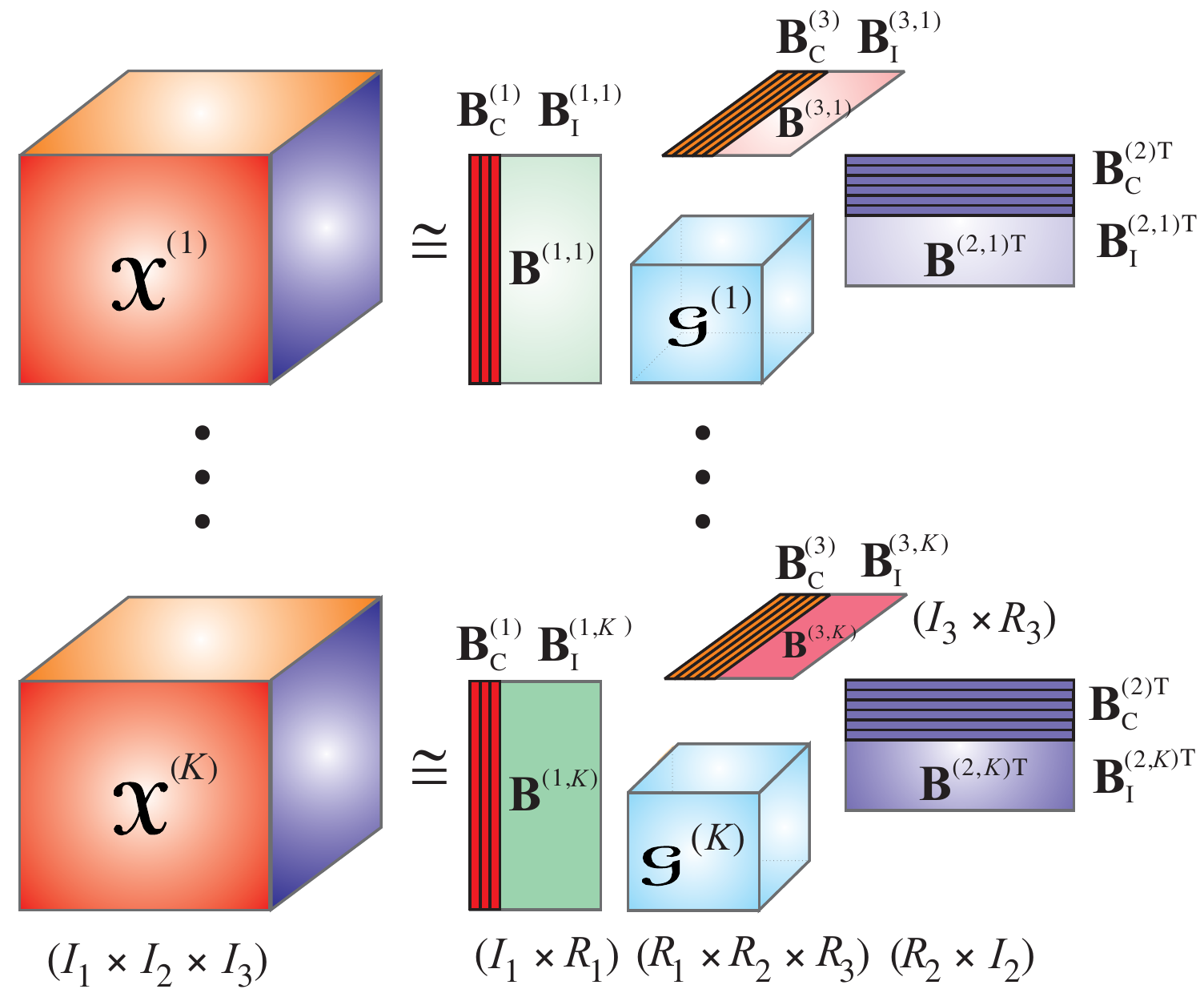}\\
 \caption{Coupled Tucker decomposition for
  linked multiway component analysis (LMWCA). The data tensors have both shared and individual components. Constraints such as orthogonality, statistical independence, sparsity and non-negativity may be imposed where appropriate. }
 \label{fig:LMBSS}
\end{figure}

The linked multiway component analysis (LMWCA) \cite{Cichocki-SICE},
shown in Figure  \ref{fig:LMBSS}, performs such decomposition into shared and individual factors, and is formulated as a set of approximate joint Tucker decompositions of a set of data tensors $ \tX^{(k)} \in  \Real^{I_1 \times I_2 \times \cdots \times I_N}$, $\;(k=1,2,\ldots,K)$:
\be
   \tX^{(k)} \cong \tG^{(k)}  \times_{1} \bB^{(1,k)}  \times_{2} \bB^{(2,k)}   \cdots  \times_{N} \bB^{(N,k)}, % + {\tE}^{(k)}, \nonumber
 \label{LMBSS}
\ee
where each factor matrix $\bB^{(n,k)}=[\bB^{(n)}_C, \; \bB^{(n,k)}_I]
\in \Real^{I_n \times R_n}$ has: (i) components $\bB^{(n)}_C \in \Real^{I_n \times C_n} $ (with $0 \leq C_n \leq R_n$) that are common (i.e., maximally correlated) to all tensors, and (ii) components $\bB^{(n,k)}_I \in \Real^{I_n \times (R_n-C_n)}$ that are tensor-specific. %The latter may, for example, correspond to processes  independent of excitations or stimuli/tasks.
The objective is to estimate the common components $\bB_C^{(n)}$, the individual components $\bB^{(n,k)}_I $, and, via the core tensors $\tG^{(k)}$, their mutual interactions. As in MWCA (see Section {\bf Tucker Decomposition}), constraints may be imposed to match data properties
 \cite{Phan2010TF,Zhou-Cichocki-MBSS}. This enables a more general and flexible framework than group ICA  and Independent Vector Analysis, which  also perform linked analysis of multiple data sets but assume that: (i) there exist only common components and (ii) the corresponding latent variables are statistically independent \cite{GICA-rev,JBSS_MCCA}, both quite stringent and limiting assumptions. As an alternative to Tucker decompositions, coupled tensor decompositions may be of a polyadic or even block term type
  \cite{Acar-all-at-once2011,Sorber-tensorlab}.

\begin{figure}[!t]
    \centering
     \includegraphics[width=10.0cm]{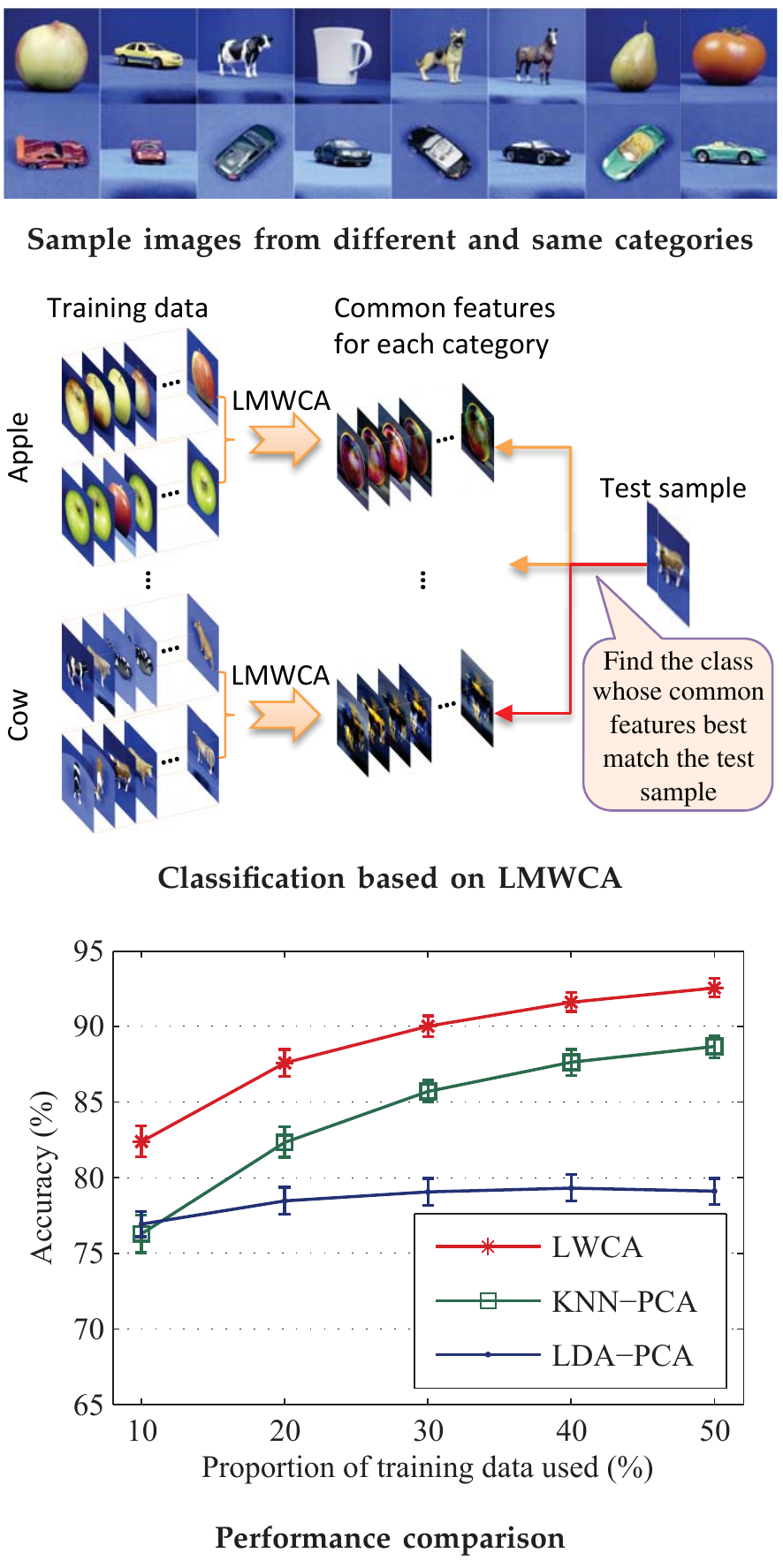}
    \caption{Classification of color objects belonging to different categories. Due to using only common features, LMWCA achieves a high classification rate, even when the training set is small.}
    \label{figLWCA}
\end{figure}

{\bf Example 5: Feature extraction and classification of objects using LMWCA.}
Classification based on common and distinct features of natural objects from the ETH-80 database  \\ (\url{http://www.d2.mpi-inf.mpg.de/Datasets}) was performed using LMWCA, whereby the discrimination among objects was performed using only the common features.
 This dataset consists of  3280 images in 8 categories, each containing 10 objects with 41 views per object.
 For each category, the training data were organized in two distinct fourth-order $(128\times128\times3\times I_4)$ tensors, where $I_4=10\times41\times0.5 p$, with  $p$ the fraction of training data.  LMWCA was applied to these two tensors to find the common and individual features, with the number of common features set to $80\%$ of $I_4$. In this way, eight sets of common features were obtained for each category. The test sample label was assigned to the category whose common features matched the new sample best (evaluated by canonical correlations) \cite{Zhou-PAMI}.
\figurename \ref{figLWCA} shows the results over 50 Monte Carlo runs and compares LMWCA with the standard  K-NN  and LDA classifiers,
the latter using  50 principal components as features. The enhanced classification results for LMWCA are attributed to the fact that the classification only makes use of the common components and is not hindered by components that are not shared across objects or views.

\section{Software}
The currently available software resources for tensor decompositions include:
\begin{itemize}
\item The Tensor Toolbox,  a versatile framework for basic operations on sparse and dense tensors, including CPD and Tucker formats \cite{tensortoolbox}.

\item The TDALAB and TENSORBOX, which provide a user-friendly interface and advanced algorithms for CPD, nonnegative Tucker decomposition and MWCA  \cite{tdalab,tensorbox}.

\item The Tensorlab toolbox builds upon the complex optimization framework and offers numerical algorithms for computing the CPD, BTD and Tucker decompositions.
The toolbox includes a library of constraints (e.g. nonnegativity, orthogonality) and the possibility to combine and jointly factorize dense, sparse and incomplete tensors \cite{Sorber-tensorlab}.

\item The $N$-Way Toolbox, which includes (constrained) CPD, Tucker decomposition and PLS in the context of chemometrics applications \cite{Nwaytoolbox}. {Many of these methods can handle constraints (e.g., nonnegativity, orthogonality) and missing elements.}

\item The TT Toolbox, the Hierarchical Tucker Toolbox and the Tensor Calculus library provide tensor tools for scientific computing \cite{oseledets2012tt,kressner2012htucker,espigtensorcalculus}.

\item Code developed for  multiway analysis is also available from the Three-Mode Company \cite{kroonenberg3modecompany}.
\end{itemize}

%\section*{Software}
%The currently available software resources for tensor decompositions include:
%\begin{itemize}
%\item The \href{<http://www.sandia.gov/~tgkolda/TensorToolbox/>}{Tensor Toolbox},  a versatile framework for basic operations on sparse and dense tensors, including CPD and Tucker formats \cite{tensortoolbox}.
%
%\item The  \href{<http://bsp.brain.riken.jp/TDALAB/>}{TDALAB} and \href{<http://www.bsp.brain.riken.jp/~phan>}{TENSORBOX}, which provide a user-friendly interface and advanced algorithms for CPD, nonnegative Tucker decomposition and MWCA  \cite{tdalab,tensorbox}.
%
%\item The \href{<http://www.esat.kuleuven.be/sista/tensorlab/>}{Tensorlab} toolbox builds upon the complex optimization framework and offers numerical algorithms for computing the CPD, BTD and Tucker decompositions.
%The toolbox includes a library of constraints (e.g. nonnegativity, orthogonality) and the possibility to combine and jointly factorize dense, sparse and incomplete tensors \cite{Sorber-tensorlab}.
%
%\item The \href{<http://www.models.life.ku.dk/nwaytoolbox>}{$N$-Way Toolbox}, which includes (constrained) CPD, Tucker decomposition and PLS in the context of chemometrics applications \cite{Nwaytoolbox}. {Many of these methods can handle constraints (e.g., nonnegativity, orthogonality) and missing elements.}
%
%\item The \href{<http://spring.inm.ras.ru/osel/?page_id=24>}{TT Toolbox}, the \href{<http://anchp.epfl.ch/htucker>}{Hierarchical Tucker Toolbox} and the \href{<http://gitorious.org/tensorcalculus>}{Tensor Calculus} library provide tensor tools for scientific computing \cite{oseledets2012tt,kressner2012htucker,espigtensorcalculus}.
%
%\item Code developed for  multiway analysis is also available from the \href{<http://three-mode.leidenuniv.nl/>}{Three-Mode Company} \cite{kroonenberg3modecompany}.
%\end{itemize}

\section{Conclusions and Future Directions}

We live in a world overwhelmed by data, from multiple pictures of Big Ben on various social web links to terabytes of data in multiview medical imaging, while \remove{in scientific experimentation} we may need to repeat the \add{scientific} experiments many times to obtain ground truth. Each snapshot gives us a somewhat  incomplete view of the same object, and involves different angles, illumination, lighting conditions, facial expressions, and noise.

We have \replace{illuminated  that}{cast a light on} tensor decompositions \replace{are}{as} a perfect match for exploratory analysis of such multifaceted data sets, and have illustrated  their applications in multi-sensor and multi-modal signal processing. Our emphasis has been to show that tensor decompositions and  multilinear algebra open completely new possibilities for component analysis, as compared with the ``flat view'' of standard two-way methods.

\replace{Contrary to}{Unlike} matrices, tensors are multiway arrays of data samples whose representations are typically overdetermined (fewer parameters in the decomposition than the number of data entries). This gives us an enormous flexibility in finding hidden components in data and the ability to enhance both robustness to noise and tolerance to missing data samples and faulty sensors. We have also discussed multilinear variants of several standard signal processing tools such as multilinear SVD, ICA, NMF and PLS, and have shown that tensor methods can operate in a deterministic way on signals of very short duration.

At present the uniqueness conditions of standard tensor models are relatively well  understood and  efficient computation algorithms do exist, however, for future applications several challenging problems remain to be addressed in more depth:
\begin{itemize}
{
\item A whole new area emerges when several decompositions which operate on different datasets are coupled, as in multiview data where some details of interest are visible in only one mode. Such techniques need theoretical support in terms of \remove{the} existence, uniqueness, and numerical properties.
\item As the complexity of advanced models increases, their computation requires efficient iterative algorithms, extending beyond the ALS class.
\item  Estimation of the number of components in data, and the assessment of their dimensionality would benefit from automation, especially in the presence of noise and outliers.
\item  Both new theory and algorithms are needed to further extend the flexibility of tensor models, e.g., for the constraints to be combined in many ways, and tailored to the particular signal properties in different modes.
\item Work on efficient techniques for saving and/or fast processing of ultra large-scale tensors is urgent, these now routinely occupy tera-bytes, and will soon require peta-bytes of memory.
\item Tools for rigorous performance analysis and rule of thumb performance bounds need to be further developed across \remove{the} tensor decomposition models.
\item Our discussion has been limited to tensor models in which all entries take values independently of one another. \remove{However, in many applications we know priors for latent variables, have prior knowledge about complex variable interaction, know that data values may belong to a certain alphabet, or can measure the noise distribution.} Probabilistic versions of tensor decompositions incorporate \remove{such} \add{prior} knowledge \add{about complex variable interaction, various data alphabets, or noise distributions}, and so promise to model data  more accurately and efficiently \cite{DinTucker,Yilmaz2010}.
  }
\end{itemize}

It is fitting to conclude  with a quote from Marcel Proust
``{\emph{The voyage of discovery is not in seeking new landscapes but in having new eyes}}''.
We  hope to have helped to bring to the eyes of the Signal Processing Community the multidisciplinary developments
in tensor decompositions, and to have shared our enthusiasm about tensors as powerful tools to discover new landscapes.
\remove{Given the current trend of different communities working on large-scale and big data analysis problems,}\add{The} future computational, visualization and interpretation tools  will be important next steps \add{in supporting the different communities working on large-scale and big data analysis problems.}

\small
%\footnotesize
\singlespacing

\end{document}